\documentclass[12pt, cite, epsf]{article}
\usepackage{amsfonts}
\usepackage{mathrsfs}

\usepackage{amstext, graphicx, amsmath, latexsym, amssymb, amscd, amsfonts}

\newtheorem{definition}{Definition}[section]
\newtheorem{theorem}{Theorem}[section]
\newtheorem{corollary}{Corollary}[section]
\newtheorem{lemma}{Lemma}[section]
\newtheorem{proposition}{Proposition}[section]
\newtheorem{remark}{Remark}[section]
\newtheorem{example}{Example}[section]

\newcommand{\R}{{\mathbb{R}}}

\DeclareMathOperator*\argmax{arg\,max}

\DeclareMathOperator*\co{conv}
\DeclareMathOperator*\cl{cl}
\DeclareMathOperator*\gph{gph}
\DeclareMathOperator*\epi{epi}

\DeclareMathOperator*\lip{lip}

\DeclareMathOperator*\rge{rge}

\begin{document}

\begin{center}
{\Large\Large\sc {\bf Kurdyka-{\L}ojasiewicz Inequality and Error Bounds of D-Gap Functions for Nonsmooth and Nonmonotone  Variational Inequality Problems}}

{\sc M.\ H.\ Li}
\\ {\small School of Mathematics and Big Data, Chongqing University of Arts and Sciences, Yongchuan, Chongqing, 402160, China }\\
Email: minghuali20021848@163.com\\[0.5cm]

{\sc K.\ W.\ Meng}
\\ {\small School of Mathematics, Southwestern University of Finance and Economics, Chengdu 611130, China}\\
Email: mengkw@swufe.edu.cn\\[0.5cm]

{\sc X. \ Q.\ Yang}\\ {\small Department of Applied Mathematics, The Hong Kong Polytechnic University, Kowloon, Hong Kong}\\
Email: mayangxq@polyu.edu.hk
\end{center}

\begin{abstract}
  In this paper, we study the D-gap function  associated with a nonsmooth and nonmonotone variational inequality problem. We present some exact formulas for the subderivative, the regular  subdifferential set, and the limiting subdifferential set of the D-gap function. By virtue of these formulas, we provide some sufficient and necessary conditions for the Kurdyka-{\L}ojasiewicz inequality property and the error bound property for the D-gap functions. As an application of our Kurdyka-{\L}ojasiewicz  inequality result and the abstract convergence result in [Attouch, et al., Convergence of descent methods for semi-algebraic and tame problems: proximal algorithms, forward-backward splitting, and regularized Gauss-Seidel methods, Math. Program., 137(2013)91-129], we show that the sequence generated by a derivative free descent algorithm with an inexact line search converges linearly to some solution of the variational inequality problem.
\end{abstract}

\noindent{\bf keywords}
Variational inequality problem,  D-gap function,  Kurdyka-{\L}ojasiewicz inequality, Error bound, Inexact line search, Linear convergence rate

\noindent{\bf AMS}
Primary, 65K10, 65K15; Secondary, 90C26, 49M37

\section{Introduction}

In this paper, we consider a  variational inequality problem (VIP) of finding $x\in K$ such that
\[
\langle F(x),y-x\rangle\geq 0 \quad\forall y\in K,
\]
where $K$ is a closed and convex subset of $\R^n$ and the mapping $F:\R^n\rightarrow \R^n$ is  locally Lipschitz continuous and  not necessarily monotone.  (VIP) has many applications in various fields such as mathematical programming, traffic network equilibrium problems and economics. We refer the reader to the very informative book \cite{fac} by Facchinei and Pang for the background information and motivations of (VIP).

 One popular approach to study (VI) is based on reformulating (VIP) as equivalent constrained/unconstrained optimization problems through the consideration of appropriate gap (merit) functions; see \cite{auc,aus,chen,fac,fuk,hark,ka,ko,lap,lig1,ng,pang,pan,pap,peng,peng1,qu,so1,tan,wu,wuzili,wuzili1,yam,yam1}. Among various reformulations in the literature, we recall that $\bar{x}$ solves (VIP) if and only if $\bar{x}$ solves the following unconstrained optimization problem with 0 as its optimal value:
\[
\min_{x\in \R^n}\quad f_{ab}(x):=f_{a}(x)-f_{b}(x),
\]
where $b>a>0$, and for each $c>0$,
\[
f_c(x):=\max_{y\in K}\{\langle F(x), x-y\rangle- \frac{c}{2}||y-x||^2\}.
\]
While $f_c$ is known as the regularized gap function \cite{auc,fuk} with  $c$ being the regularized parameter, $f_{ab}$ is often known as the D-gap function \cite{peng} with `D' standing for the 'difference' of two parameterized regularized gap functions.  By replacing the quadratic term in defining $f_c$ with some general term having very similar properties as those of the quadratic term, the corresponding {\it generalized} regularized gap  and  {\it generalized} D-gap functions have also been extensively studied in the literature; see \cite{lig,lig1,wu,yam1}.

The (generalized) differentiability properties of these regularized gap and D-gap functions have been extensively investigated, and have been utilized to study the property of  error bounds \cite{fac} and the property of the  Kurdyka-{\L}ojasiewicz (KL, for short) inequality  \cite{dp22}. The latter properties have played very important roles in convergence analysis for algorithms designed based upon gap functions.

We review a few typical  results  related to the (generalized) D-gap function as follows. Peng \cite{peng} showed  that if $F$ is  continuously differentiable and strongly monotone,    the D-gap function is also continuously differentiable and its square root provides a global error bound for (VIP).
Yamashita et al. \cite{yam1} introduced the generalized D-gap function and obtained its continuous differentiability by assuming that $F$ is continuously differentiable.  Moreover, by assuming that $F$ is strongly monotone and that either $F$ is Lipschitz continuous or $K$ is compact, they showed that the square root of the  generalized D-gap function provides a global error bound for (VIP), and that the sequence generated by a descent algorithm with an inexact line search converges  to the unique solution of (VIP).
 Based on the D-gap function and by assuming that $F$ is continuously differentiable and monotone, Solodov and Tseng \cite{so} developed two unconstrained methods that are similar to the feasible method in Zhu and Marcotte \cite{zhu} which is based on the regularized gap function.
 By assuming that $F$ is locally Lipschitz continuous,  Xu \cite{xu} obtained a formula for  the Clarke subdifferential set of the D-gap function,  and a   global convergence result for a descent algorithm with an inexact line search under the circumstance that $F$ is strongly monotone and Lipschitz continuous.
     By the same assumption that  $F$ is locally Lipschitz continuous, Ng and Tan \cite{ng2} obtained some formulas for the Clarke directional derivative and the Clarke subdifferential set of the D-gap function.
 By assuming that $F$ is coercive and locally Lipschitz continuous, and by introducing a condition expressed in terms of  the Clarke generalized Jacobian of $F$, Li and Ng \cite{lig} showed that the square root of the generalized D-gap function provides a local error bound for (VIP),  and by virtue of which, they proved that any cluster point of the sequence generated by a descent algorithm with an inexact line search is a solution of (VIP), and that the convergence rate is linear when $F$ is smooth,  strongly monotone and $\nabla F$ is locally Lipschitz continuous. Note that Li and Ng \cite{lig} also provided some formulas for the Clarke directional derivative and the Clarke subdifferential set of the generalized D-gap function, which were very  crucial for their arguments.
   Later Li et al. \cite{lig1} established some  error bound results for the generalized D-gap function  by assuming that  $F$ is (Lipschitz) continuous, locally monotone and coercive.

From the literature review above, it is clear to see that most of the existing results for error bounds and the convergence of a descent algorithm were obtained by assuming that $F$ is strongly monotone, with an exception being that,  the error bound result in Li and Ng \cite{lig}, though having difficulty in verification, was applied to some cases when $F$ is nonmonotone. As for the property of the KL inequality,  there is almost no result, to the best of our knowledge,  presented in a straightforward way for the case when  $F$ is  locally Lipschitz continuous.
By examining the definition for the KL inequality (see Definition \ref{def-kl}  below) and the theory of error bounds in \cite{ac04, my12}, it is reasonable that the notion of the subderivative, the regular/Fr\'{e}chet  subdifferential set,  and the  general/limiting subdifferential set (see Definition \ref{def-guangyi-daoshu}) should have  played a role in studying the generalized differentiability properties of the regularized gap and D-gap functions. But it is quite surprising that there is no such a related result in the literature for the case when $F$ is locally Lipschitz continuous and not necessarily monotone.

To fill this gap, we will investigate the KL inequality and error bounds of the D-gap function for nonsmooth and nonmonotone (VIP) by providing formulas for the subderivative and the (limiting) subdifferential sets of the D-gap functions, and as an application of our result for the KL inequality and the abstract convergence result in \cite{bot1} for inexact descent methods, we will establish the linear convergence rate for a descent algorithm with an inexact line search.

The main contributions of the  paper are  as follows.
\begin{description}
  \item[(i)] We obtain a number of exact formulas for the subderivatives, the regular/Fr\'{e}chet  subdifferential sets, and the general/limiting subdifferential sets of the regularized gap function $f_c$ and the D-gap function $f_{ab}$, respectively. See Propositions \ref{subderivative000}-\ref{prop-sub-diff} below.
      Taking the limiting subdifferential set $\partial f_{ab}(\bar{x})$ of $f_{ab}$ at a point $\bar{x}$ for instance, we obtain
      $$
      \partial f_{ab}(\bar{x})=\displaystyle D^{*}F(\bar{x}) \left(\pi_b(\bar{x})-\pi_a(\bar{x})\right)-b\left(\pi_b(\bar{x})-\pi_a(\bar{x})\right)+(b-a)
(\bar{x}-\pi_a(\bar{x})),
$$
where $D^{*}F(\bar{x})$ denotes the coderivative of $F$ at $\bar{x}$ (cf. Definition \ref{def-set-valued-derivative}), and for each $c>0$,  $\pi_{c}(x):=P_K\left(x-\frac{F\left(x\right)}{c}\right)$ with $P_K(\cdot)$ being the projection operator onto $K$.
  To the best of our knowledge, these formulas have not been seen from the literature, although, as mentioned above, 
  exact formulas have been obtained for the Clarke directional derivatives and the Clarke subdifferential sets of
   $f_c$ and $f_{ab}$, respectively. It should be noticed that, although $f_c$ is a marginal function and $f_{ab}=f_a-f_b$ is a difference of two marginal functions, we cannot obtain our formulas by directly applying the theory of marginal
    functions known from the literature \cite{aaus,hogan,minsak,morny,rocka}.  As a matter of fact,
     our approach depends heavily on the inherent structures of $f_c$ and $f_{ab}$.   
  \item[(ii)] By virtue of the formula obtained for the general/limiting subdifferential set  of the D-gap function $f_{ab}$, we present a few sharp results on the properties of the KL inequality and the error bounds for $f_{ab}$. In particular, by assuming that the following inequality holds for some $\mu>0$ and for all
  $x\in \R^n$ where $F$ is differentiable:
 \begin{equation}\label{key-suff}
 \displaystyle\langle \nabla F(x) (\pi_a(x)-\pi_b(x)),\;\pi_a(x)-\pi_b(x)\rangle\geq \mu||\pi_a(x)-\pi_b(x)||^2,
\end{equation}
 which can be considered as a restricted (weaker) notion of strong monotonicity, we show that
 $$
 d(0, \partial f_{ab}(x) )\geq \mu \|\pi_b(x)-\pi_a(x)\| \quad \forall x\in \R^n,
 $$
 and  that $f_{ab}$ is a KL  function  with an exponent of $\frac{1}{2}$,  and moreover  that some local/global error bound results holds. See Theorem \ref{theo-zuizhongyao} below.
  \item[(iii)] By assuming (\ref{key-suff}) and applying our result on the KL property for $f_{ab}$, we obtained  the linear convergence rate for a derivative free descent algorithm, which is essentially the same algorithm as those studied in \cite{ka,lig,peng1,qu,xu,yam1}. See Theorem \ref{convergence} below.
      Starting from any initial point $x_0$, the algorithm generates a sequence $\{x_n\}$ in the manner of $x_{n+1}=x_{n}+t_nd_n$, where $d_n$ is the search direction, either being $\pi_{a}(x_n)-x_n$ or $\pi_{a}(x_n)-\pi_{b}(x_n)$, and $t_n$ is the stepsize determined by an Armijo line search. Under some other mild assumptions, except for (\ref{key-suff}), we show that the stepsize sequence $\{t_n\}$ has a positive lower bound $t^*>0$ (cf. Proposition \ref{positive bound} below), and moreover the following hold (cf. Proposition \ref{prop-key conditions} below):
\[
f_{ab}(x_{n+1})-f_{ab}(x_n)\leq  -M_1||x_{n+1}-x_n||^2
\]
 and
 \[
d(0, \partial f_{ab}(x_n))\leq    \frac{M_2}{t^*} ||x_{n+1}-x_n||,
\]
where $M_1$ and $M_2 $ are two positive constants.  That is,  the sequence $\{x_n\}$  satisfies the assumptions {\bf (H1)}, a variant of {\bf(H2)},
 and {\bf (H3)} proposed in \cite{bot1},  and our convergence analysis falls into the framework of the abstract convergence  for inexact descent methods studied in \cite{bot1}.
\end{description}

The outline of the paper is as follows. Section \ref{sec-mp} is about notation and terminology, and some mathematical preliminaries.  In section \ref{sec-ss-d-gap}, we present some exact formulas for the subderivatives, the regular/Fr\'{e}chet  subdifferential sets, and the general/limiting subdifferential sets of the regularized gap function $f_c$ and the D-gap function $f_{ab}$, respectively. By virtue of these formulas for the D-gap function, we present in Section \ref{sec-kl-error} some sufficient and necessary conditions for the error bound property and the KL inequality property. As an application of our KL inequality result  and the abstract convergence result in \cite{bot1} for inexact descent methods, we show in section \ref{sec-method} that the sequence generated by a descent algorithm (based upon the D-gap function) with an inexact line search converges linearly to some solution of (VIP).

\section{Notation and Mathematical Preliminaries}\label{sec-mp}
Throughout the paper we use the standard notations of variational analysis; see the seminal book \cite{roc} by Rockafellar and Wets.
 The Euclidean norm of a vector $x$ is denoted by $||x||$, and the inner product of vectors $x$ and $y$ is denoted by $\langle x, y\rangle$.
  Let $A\subset \R^n$ be a nonempty set.
We denote   by $\co A$ the convex hull of $A$. The polar cone of $A$ is defined by
$
A^*:=\{v\in \R^n\mid \langle v, x\rangle\leq 0\;\forall x\in A\}.
$
The distance from $x$ to $A$ is defined by
$
d(x,A):=\inf_{y\in A}||y-x||.
$
The projection mapping $P_A$ is defined by
$
P_A(x):=\{y\in A\mid \|y-x\|=d(x,A)\}.
$

 \begin{definition}\label{def-jihe-C}
  Let $C\subset \R^n$  and let $x\in C$.
 \begin{description}
   \item[{\bf (i)}]  The  tangent cone to $C$ at $x$ is denoted by $T_C(x)$, i.e., $w\in T_C(x)$ if there exist sequences $t_k\downarrow 0$ and $\{w_k\}\subset \R^n$ with $w_k\rightarrow w$ and $x+ t_k w_k\in C \;\forall k$.
   \item[{\bf (ii)}] The regular   normal cone to $C$ at $x$ is denoted by $\widehat{N}_C(x)$, i.e.,  $v\in \widehat{N}_C(x)$ if
   \[
   \langle v, x-\bar{x}\rangle \leq o(\|x-\bar{x}\|)\quad\mbox{for all } x\in C.
   \]
   Another way of defining the regular normal cone  is via the equality $\widehat{N}_C(x)=T_C(x)^*$.
   \item[{\bf (iii)}] The   normal cone  to $C$ at $x$ is denoted by $N_C(x)$, i.e., $v\in N_C(x)$ if  there exist sequences $x_k\to x$ and $v_k\to v$ with $x_k\in C$ and $v_k\in \widehat{N}_C(x_k)$ for all $k$.
 \item[{\bf (iv)}] $C$ is said to be regular at $x$ in the sense of Clarke if it is locally closed at $x$ (i.e., $C\cap U$ is closed for some closed neighborhood $U$ of $x$) and $\widehat{N}_C(x)=N_C(x)$.
 \end{description}
 \end{definition}

Let $f:\R^n\to \overline{\R}:=\R\cup\{\pm \infty\}$ be an extended real-valued function.  We denote the epigraph of $f$ by
$
\epi f:=\{(x,\alpha)\mid f(x)\leq \alpha\}.
$
The lower level set with a  level of $\alpha$ is defined and denoted by
$
[f\leq \alpha]:=\{x\in \R^n\mid f(x)\leq \alpha\}.
$
In a similar way, we define
$[f< \alpha]:=\{x\in \R^n\mid f(x)< \alpha\}$  and $[\alpha<f< \beta]:=\{x\in \R^n\mid \alpha<f(x)< \beta\}.$

\begin{definition}\label{def-guangyi-daoshu}
Let $f:\R^n\to \overline{\R}$ be an extended real-valued function and let $\bar{x}$ be  a point with $f(\bar{x})$ finite.
\begin{description}
  \item[{\bf (i)}] The vector $v\in \R^n$ is a regular/Fr\'{e}chet  subgradient of $f$ at $\bar{x}$, written $v\in \widehat{\partial} f(\bar{x})$, if
\[
f(x)\geq f(\bar{x})+\langle v, x-\bar{x}\rangle+o(||x-\bar{x}||).
\]
  \item[{\bf (ii)}] The vector $v\in \R^n$ is a  general/limiting subgradient of $f$ at $\bar{x}$, written $v\in \partial f(\bar{x})$, if there exist sequences $x_k\to \bar{x}$ and $v_k\to v$ with $f(x_k)\to f(\bar{x})$ and $v_k\in \widehat{\partial} f(x_k)$.
    \item[{\bf (iii)}] The function $f$ is said to be (subdifferentially) regular at $\bar{x}$ if $\epi f$ is regular
 in the sense of Clarke at $(\bar{x}, f(\bar{x}))$ as a subset of $\R^n\times \R$.
 \item[{\bf (iv)}]  The subderivative $df(\bar{x}):\R^n\to \overline{\R}$ is defined by
$$
df(\bar{x})(w):=\liminf_{t\downarrow 0, w'\to w}\frac{f(\bar{x}+tw')-f(\bar{x})}{t}.
$$
\item[{\bf (v)}] The set of Clarke subgradients of $f$ at $\bar{x}$ is defined by
\[
\overline{\partial} f(\bar{x}):=\{v| (v,-1)\in \cl \co N_{\epi f}(\bar{x},f(\bar{x}))\},
\]
where $\cl \co N_{\epi f}(\bar{x},f(\bar{x}))$ denotes the closed and convex hull of $N_{\epi f}(\bar{x},f(\bar{x}))$.
\end{description}
\end{definition}

\begin{remark}\label{rem-reg-sub}
The regular subgradients can be derived from the subderivative as follows \cite[Exercise 8.4]{roc}:
$$
\widehat{\partial} f(\bar{x})=\{v\in \R^n| \langle v, w\rangle\leq df(\bar{x})(w)\;\forall w\in \R^n\}.
$$
\end{remark}

Following \cite{bot0,bot2,lipong}, we introduce the notion of the Kurdyka-{\L}ojasiewicz (KL, for short) inequality.
\begin{definition}\label{def-kl}
For a proper lower semicontinuous  function $f:\R^n\to \overline{\R}:=\R\cup \{\pm \infty\}$, a point $\bar{x}\in \R^n$ with $\partial f(\bar{x})\not=\emptyset$,   and some $\alpha\in [0, 1)$,  we say that  $f$ satisfies the KL inequality   at $\bar{x}$ with an exponent of $\alpha$, if  there exist   $\mu, \epsilon>0$ and  $\nu\in (0,+\infty]$ so that
\[
d(0, \partial f(x))\geq \mu(f(x)-f(\bar{x}))^\alpha
\]
whenever $\|x-\bar{x}\|\leq \epsilon$ and $f(\bar{x})<f(x)<f(\bar{x})+\nu$. If f satisfies the KL inequality  at every $x\in \R^n$ with $\partial f(x)\not=\emptyset$ and with   the same exponent $\alpha$,  we say that $f$ is a KL function with an exponent of $\alpha$.
\end{definition}

Following \cite{fac}, we introduce the notion of local and global error bounds as follows.
\begin{definition}\label{def-error-bound}
For a proper function $f:\R^n\to \overline{\R}$ and a set $C\subset \R^n$,   we say that  $f$  has a local error bound  on $C$  if  there exist two positive constants   $\tau$ and  $\epsilon$ such that for all $x\in [f\leq \epsilon]\cap C$
\[
d(x,  [f\leq 0]\cap C)\leq \tau \max\{f(x), 0\}.
\]
Furthermore,  we say that $f$ has a global error bound on $C$ if there exists a constant $\tau>0$ such that the above inequality holds for all $x\in C$.
\end{definition}

  \begin{definition}\label{def-set-valued-derivative}
  Let $S:\R^n\rightrightarrows \R^m$ be a set-valued mapping and $(\bar{x}, \bar{u})\in \gph S:=\{(x, u)\mid u\in S(x)\}.$
  \begin{description}
    \item[{\bf (i)}]  The graphical   derivative of $S$ at $\bar{x}$ for $\bar{u}$ is the mapping $DS(\bar{x}\mid \bar{u}):\R^n\rightrightarrows \R^m$ defined by
     \[
z\in DS(\bar{x}\mid \bar{u})(w)\Longleftrightarrow (w, z)\in T_{\gph S}(\bar{x}, \bar{u}).
  \]
 \item[{\bf (ii)}] The regular coderivative of $S$ at $\bar{x}$ for $\bar{u}$  is the mapping $\widehat{D}^*S(\bar{x}\mid \bar{u}):\R^m\rightrightarrows \R^n$ defined by
  \[
  x^*\in \widehat{D}^*S(\bar{x}\mid \bar{u})(u^*)\Longleftrightarrow (x^*, -u^*)\in \widehat{N}_{\gph S}(\bar{x}, \bar{u}).
  \]
     \item[{\bf (iii)}]  The    coderivative of $S$ at $\bar{x}$ for $\bar{u}$ is
    the mapping $D^*S(\bar{x}\mid \bar{u}):\R^m\rightrightarrows \R^n$ defined by
  \[
  x^*\in D^*S(\bar{x}\mid \bar{u})(u^*)\Longleftrightarrow (x^*, -u^*)\in N_{\gph S}(\bar{x}, \bar{u}).
  \]
  \end{description}
    Here the notation $DS(\bar{x}\mid \bar{u})$, $D^*S(\bar{x}\mid \bar{u})$ and $\widehat{D}^*S(\bar{x}\mid \bar{u})$ is simplified to $DS(\bar{x})$, $D^*S(\bar{x})$ and $\widehat{D}^*S(\bar{x})$ when $S$ is single-valued at $\bar{x}$, i.e., $S(\bar{x})=\{\bar{u}\}$.
  \end{definition}

\begin{definition}\label{def-lip-F}
Let $F$ be a single-valued mapping defined on $\R^n$, with values in $\R^m$.
\begin{description}
  \item[{\bf (i)}] $F$ is  globally  Lipschitz continuous  if there exists $\kappa \in \R_+:=[0, \infty)$ with
  \[
  \|F(x')-F(x)\|\leq \kappa\|x'-x\|\quad \forall x,x'\in \R^n.
  \]
  Then $\kappa$ is called a Lipschitz constant for $F$.
  \item[{\bf (ii)}] $F$ is locally Lipschitz continuous at a point  $\bar{x}\in \R^n$ if   the value
   \[
  \lip F(\bar{x}):=\limsup_{x,x'\to \bar{x}, x\not=x'}\frac{\|F(x')-F(x)\|}{\|x'-x\|}
  \]
  is finite. Here  ${\lip}\, F(\bar{x})$ is the Lipschitz modulus of $F$ at $\bar{x}$.
\item[{\bf (iii)}] $F$ is locally Lipschitz continuous if  $F$ is locally Lipschitz continuous at every $\bar{x}\in \R^n$.
 \end{description}
 \end{definition}


\begin{lemma}\label{lem-lip-f-pro}
Let $f:\R^n\to \overline{\R}$ be an extended real-valued function and let $\bar{x}$ be  a point with $f(\bar{x})$ finite. Assume that $f$ is locally Lipschitz continuous at $\bar{x}$. The following properties hold:
\begin{description}
  \item[{\bf (a)}] $\partial f(\bar{x})$ is nonempty and compact.
  \item[{\bf (b)}] $df(\bar{x})(w)=\displaystyle\liminf_{t\downarrow 0}\frac{f(\bar{x}+tw)-f(\bar{x})}{t}$.
  \item[{\bf (c)}] $\overline{\partial} f(\bar{x})=\co(\partial f(\bar{x}))$.
\end{description}
\end{lemma}
\noindent{\bf Proof.} {\bf (a-c)} can be found in \cite[ Theorem 9.13, Exercise 9.15, Theorem 9.61]{roc}, respectively.
\hfill $\Box$

\begin{lemma}\label{lem-lip-F-pro}
Assume that $F:\R^n\to \R^m$ is locally Lipschitz continuous at a point $\bar{x}\in \R^n$. The following properties hold:
\begin{description}
  \item[{\bf (a)}] $D^*F(\bar{x})(0)=\{0\}$, which is also sufficient for $F$ being locally Lipschitz continuous at $\bar{x}$.
  \item[{\bf (b)}]  The mappings  $DF(\bar{x})$ and $D^*F(\bar{x})$ are nonempty-valued and locally bounded. 
  \item[{\bf (c)}] $||z||\leq (\lip F(\bar{x}))\, ||w||$ holds for all $(w,z)\in \gph (DF(\bar{x}))$.
   \item[{\bf (d)}] $||x^*||\leq (\lip F(\bar{x}))\, ||u^*||$ holds for all $(u^*, x^*)\in \gph(D^*F(\bar{x}))$. 
\item[{\bf (e)}] $z\in DF(\bar{x})(w)$ if and only if there is some $\tau^\nu\downarrow 0$ such that
$
 \frac{F(\bar{x}+\tau^\nu w)-F(\bar{x})}{\tau^\nu}\rightarrow z.
$
\end{description}
\end{lemma}
\noindent{\bf Proof.} {\bf (a)} follows directly from the Mordukhovich criterion \cite[Theorem 9.40]{roc}. {\bf (b-d)} follow  from \cite[Proposition 9.24]{roc}. {\bf (e)} follows from the definitions of the graphical derivative and the local Lipschitzian continuity.  \hfill$\Box$

Assume now that $F:\R^n\to \R^m$ is a locally Lipschitz continuous function and let $D$ be the subset of $\R^n$ consisting of the points where $F$ is differentiable. By the Rademacher Theorem \cite[Theorem 9.60]{roc}, $F$ is differentiable almost everywhere with  $\R^n\backslash D$ being negligible. For each $\bar{x}\in \R^n$, define
\begin{equation}\label{def-B-jac}
\overline{\nabla} F(\bar{x}):=\{A\in \R^{m\times n}\mid \exists x^\nu\to \bar{x}\;\mbox{with}\;x^\nu\in D,\,\nabla F(x^\nu)\to A\},
\end{equation}
in terms of which,  the generalized Jacobian $\overline{\partial} F(x)$   \cite[Definition 2.6.1]{cla} of $F$ at $\bar{x}$ can be written as
\begin{equation}\label{g-jac}
\overline{\partial} F(\bar{x}):=\co  \overline{\nabla} F(\bar{x}).
\end{equation}
According to \cite[Theorem 9.62]{roc},  $\overline{\nabla} F(\bar{x})$ is a nonempty, compact set of matrices, and for every $w\in \R^n$ and $y\in \R^m$ one has
\begin{equation}\label{code-dual}
\co D^*F(\bar{x})(y)=\co\{A^Ty\mid A\in \overline{\nabla} F(\bar{x})\}=\{A^Ty\mid A\in \co \overline{\nabla} F(\bar{x})\}
\end{equation}
and
\begin{equation}\label{stri-dual}
\co D_*F(\bar{x})(w)=\co\{Aw\mid A\in \overline{\nabla} F(\bar{x})\}=\{Aw\mid A\in \co \overline{\nabla} F(\bar{x})\},
\end{equation}
 where $D_*F(\bar{x})$ stands for the strict derivative mapping of $F$ at $\bar{x}$ \cite[Definition 9.53]{roc}, and has the following definition by taking into account that $F$ is locally Lipschitz continuous:
\begin{equation}\label{def-strict-der}
D_*F(\bar{x})(w):=\{z\mid \exists \tau^\nu \downarrow 0, x^\nu\to \bar{x}\;\mbox{with}\;(F(x^\nu+\tau^\nu w)-F(x^\nu))/\tau^\nu\to z\}.
\end{equation}
Note that $D_*F(\bar{x})$ is also known as the Thibault's strict derivative (cf. \cite{th}), and that by definition
\begin{equation}\label{txds-strict}
\gph DF(\bar{x})\subset \gph D_*F(\bar{x}).
\end{equation}

\begin{definition} \cite{fac}\label{def-mono-coer}
Let $C$ be a subset of $\R^n$, and let $F$ be a single-valued mapping defined on $\R^n$, with values in $\R^n$. $F$ is said to be coercive on $C$ if
$$\lim_{x\in C,\,\|x\|\to \infty}\frac{\langle F(x), x-y\rangle}{\|x\|}=+\infty$$ holds for all $y\in C$
(if $C$ is bounded, then  $F$ is by convention coercive on $C$);  and $F$ is said to be
    strongly monotone on $C$ (with modulus $\mu>0$) if
  $
  \langle F(x)-F(y),\; x-y\rangle \geq \mu\|x-y\|^2$ holds for all $x, y\in C$.
 \end{definition}

\section{Subderivatives and subgradients of gap functions}\label{sec-ss-d-gap}

In the remainder of the paper,  we make  the following blanket assumptions on problem data and some constants,  and for the sake of simplicity,  we will not mention them in stating a result.
 \begin{itemize}
   \item $K\subset \R^n$ is a nonempty closed and convex set.
   \item $F:\R^n\to \R^n$ is a locally Lipschitz continuous function.
   \item $ a, b, c$  are fixed positive numbers with $a<b$.
 \end{itemize}

The aim of this section is to study subderivatives and subgradients of $f_{ab}$ and $f_c$ at some $\bar{x}$
 by virtue of the graphical derivative $DF(\bar{x})$  and the coderivatives,  $D^*F(\bar{x})$ and $\widehat{D}^*F(\bar{x})$, and frequently,
the following projection operator $\pi_c$ associated with $F$ and $K$:
$$
\pi_{c}(x):=P_K\left(x-\frac{F\left(x\right)}{c}\right).
$$
The projection operators $\pi_a$ and $\pi_b$ are defined in the same manner.  

To begin with, we summarize below  some basic properties of  the regularized gap function $f_c$ and the D-gap function $f_{ab}$,
most of which can be found in the literature and are useful for further development in the sequel.
  \begin{lemma}\label{lem-basic}
 The following  properties  hold:
\begin{description}
\item[{\bf (a)}]$\frac{b-a}{2}||x-\pi_{b}(x)||^2+\frac{a}{2}||\pi_{b}(x)-\pi_{a}(x)||^2\leq f_{ab}(x)\leq \frac{b-a}{2}||x-\pi_{a}(x)||^2-\frac{b}{2}||\pi_{b}(x)-\pi_{a}(x)||^2$.
\item[{\bf (b)}] $||\pi_{b}(x)-\pi_{a}(x)||\leq \frac{b-a}{a}||x-\pi_{a}(x)||$ and $||x-\pi_{b}(x)||\leq ||x-\pi_{a}(x)||\leq \frac{b}{a}||x-\pi_{b}(x)||$.
\item[{\bf (c)}] $x\in \R^n$ solves (VIP) $\Leftrightarrow$ $x=\pi_c(x)$ for any $c>0$  $\Leftrightarrow$   $f_{ab}(y)\geq f_{ab}(x)=0$ for all $y\in \R^n$ $\Leftrightarrow$   $x\in K$ and $f_c(y)\geq f_c(x)=0$ for all $y\in K$.
\item[{\bf (d)}]  $\langle a(x-\pi_a(x))-b(x-\pi_b(x)), \pi_a(x)-\pi_b(x)\rangle\geq 0$.
\item[{\bf (e)}] $\pi_a(x)-\pi_b(x)\in T_{ab}(x,F,K):=T_K(\pi_b(x))\cap (-T_K (\pi_a(x)))\cap (F(x))^{*}$.
\item[{\bf (f)}]  $\pi_a$, $\pi_b$, $\pi_c$,  $f_c$ and $f_{ab}$  are  locally Lipschitz continuous. If $F$ is globally Lipschitz continuous, then $\pi_a$, $\pi_b$, $\pi_c$,  $f_c$ and $f_{ab}$ are also  globally Lipschitz continuous.
\item[{\bf (g)}] The following hold:
\[
\begin{array}{ll}
&\argmax_{y\in K}\left\{\left\langle F(x), x-y\right\rangle- \frac{c}{2}||y-x||^2\right\}=\{\pi_{c}(x)\},\\[0.1cm]
&f_{c}(x)=\langle F(x), x-\pi_{c}(x)\rangle-\frac{c}{2} ||x-\pi_{c}(x)||^2,\\[0.1cm]
&f_{ab}(x)=\langle F(x), \pi_{b}(x)-\pi_{a}(x)\rangle-\frac{a}{2} ||x-\pi_{a}(x)||^2+\frac{b}{2} ||x-\pi_{b}(x)||^2.
\end{array}
\]
\end{description}
\end{lemma}
\noindent{\bf Proof.}  {\bf (a)} and {\bf (b)} can be found in \cite[Lemma 1]{so} and \cite{ng2}, respectively.
 {\bf (c)} can be found in \cite{fuk} and \cite{wu}.   {\bf (d)} and {\bf (e)} can be found in \cite[Lemma 4.4]{lig} or in \cite[Theorem 10.3.4]{fac}.
  {\bf (f)} can be found in \cite[Lemma 3.1]{lig1}.    {\bf (g)} can be found in  \cite{wu} or deduced from standard optimality condition for convex  programs.  This completes the proof. \hfill $\Box$

\subsection{Subderivatives and subgradients of $f_c$}\label{sec-gamma}
We first present the formulas for the subderivative, the regular   subdifferential set and the limiting    subdifferential set of $f_c$ at a point $\bar{x}$.
\begin{proposition}\label{subderivative000}
Let $\bar{x}\in \R^n$ and let $w\in \R^n$. We have the following formulas:
\[
\begin{array}{ll}
&df_{c}(\bar{x})(w)=\displaystyle\langle F(\bar{x}),\;w\rangle +\min \langle \left(DF(\bar{x})-c I\right) w, \; \bar{x}-\pi_c(\bar{x})\rangle,\\[0.1cm]
&\widehat{\partial}f_{c}(\bar{x})=\left(\widehat{D}^{*}F(\bar{x})-c I\right)\left( \bar{x}-\pi_c(\bar{x}) \right)+F(\bar{x}),\\[0.1cm]
&\partial f_c(\bar{x})=\displaystyle \left(D^{*}F(\bar{x})-c I\right)\left( \bar{x}-\pi_c(\bar{x}) \right)+F(\bar{x}),
\end{array}
\]
where $$\min \langle \left( DF(\bar{x})-cI \right)w, \; \bar{x}-\pi_c(\bar{x})\rangle:=\displaystyle\min_{v\in DF(\bar{x})(w)} \langle v-cw, \; \bar{x}-\pi_c(\bar{x})\rangle.$$
\end{proposition}
\noindent{\bf Proof.}  Let $w\in \R^n$ be fixed.   Since $F$ is locally Lipschitz continuous,  it follows from Lemma \ref{lem-lip-F-pro} {\bf (b)} and {\bf (e)} that for any continuous function $M:\R\to \R^n$,
\begin{equation}\label{chouxiang}
\liminf_{t\downarrow 0}\langle \frac{F(\bar{x}+tw)-F(\bar{x})}{t},\; M(t)\rangle=\min_{v\in DF(\bar{x})(w)}\langle v,\; M(0)\rangle.
\end{equation}
By Lemma \ref{lem-basic} {\bf (f)}, $f_c$ is a locally Lipschitz continuous function, which implies by Lemma \ref{lem-lip-f-pro} {\bf (b)} that
$
df_c(\bar{x})(w)=\liminf_{t\downarrow 0} \frac{f_c(\bar{x}+t w)-f_c(\bar{x})}{t}$.
In view of Lemma \ref{lem-basic} {\bf (g)}, we have for all $t$,
$f_c(\bar{x})\geq \langle F(\bar{x}),\bar{x}-\pi_c(\bar{x}+t w)\rangle-\frac{c}{2} ||\bar{x}-\pi_c(\bar{x}+tw)||^2$, and
$f_c(\bar{x}+t w)= \langle F(\bar{x}+tw),\bar{x}+tw-\pi_c(\bar{x}+t w)\rangle-\frac{c}{2} ||\bar{x}+tw-\pi_c(\bar{x}+t w)||^2$.
This, together with (\ref{chouxiang}) and the fact that $\pi_c$ is locally Lipschitz continuous (cf. Lemma \ref{lem-basic} {\bf (f)}), implies that
\begin{eqnarray}
df_{c}(\bar{x})(w)&\leq& \liminf_{t\downarrow 0}\langle \frac{F(\bar{x}+tw)-F(\bar{x})}{t},\; \bar{x}-\pi_c(\bar{x}+tw)\rangle+\lim_{t\downarrow 0}\langle F(\bar{x}+tw),\; w\rangle\nonumber\\[0.1cm]
&&+\lim_{t\downarrow 0}\frac{c}{2}\langle 2(\bar{x}-\pi_c(\bar{x}+tw))+tw,\;-w\rangle\nonumber\\[0.1cm]
&=& \min_{v\in DF(\bar{x})(w)}\langle v,  \bar{x} -\pi_c(\bar{x})\rangle+\langle F(\bar{x}),\; w\rangle-c\langle  \bar{x} -\pi_c(\bar{x}),\; w\rangle\nonumber\\[0.1cm]
&=& \min \langle \left(DF(\bar{x})-c I\right) w, \; \bar{x}-\pi_c(\bar{x})\rangle+\displaystyle\langle F(\bar{x}),\;w\rangle.\nonumber
\end{eqnarray}
To prove the inequality in the other direction, we simply follow a similar way by observing  from Lemma \ref{lem-basic} {\bf (g)} that for all $t$, $
f_c(\bar{x})= \langle F(\bar{x}),\bar{x}-\pi_c(\bar{x})\rangle-\frac{c}{2} ||\bar{x}-\pi_c(\bar{x})||^2$, and
$f_c(\bar{x}+t w)\geq \langle F(\bar{x}+tw),\bar{x}+tw-\pi_c(\bar{x})\rangle-\frac{c}{2} ||\bar{x}+tw-\pi_c(\bar{x})||^2$.

To get the formula for $\widehat{\partial}f_{c}(\bar{x})$, we resort to the formula for $df_{c}(\bar{x})$ and  the equality
in Remark \ref{rem-reg-sub}.
Specifically, in terms of  $\bar{v}:=F(\bar{x})-c (\bar{x}-\pi_c(\bar{x}))$,  we have
\[
\begin{array}{lll}
&&v\in \widehat{\partial} f_{c}(\bar{x})\\[0.1cm]
&\Longleftrightarrow&\langle v, w\rangle\leq \langle \bar{v}, w\rangle+\min \langle DF(\bar{x})(w),  \bar{x} -\pi_c(\bar{x})\rangle\quad\forall w\in \R^n,\\[0.1cm]
&\Longleftrightarrow& \langle v-\bar{v}, w\rangle\leq \langle z, \bar{x} -\pi_c(\bar{x})\rangle\quad\forall (w, z)\in \gph (DF(\bar{x}))= T_{ \gph F} (\bar{x}, F(\bar{x})),\\[0.1cm]
&\Longleftrightarrow& (v-\bar{v}, -\bar{x}+\pi_c(\bar{x}))\in (T_{ \gph F} (\bar{x}, F(\bar{x})))^* =\widehat{N}_{\mbox{gph} F}(\bar{x}, F(\bar{x})),\\[0.1cm]
&\Longleftrightarrow& v-\bar{v}\in \widehat{D}^{*}F(\bar{x})(\bar{x} -\pi_c(\bar{x})).
\end{array}
\]
This gives us the formula for $\widehat{\partial}f_{c}(\bar{x})$.

 To show  $\partial f_{c}(\bar{x})\subset U:=\displaystyle \left(D^{*}F(\bar{x})-c I\right)\left( \bar{x}-\pi_c(\bar{x}) \right)+F(\bar{x})$,  let $v\in \partial f_{c}(\bar{x})$.  Then by the formula for $\widehat{\partial}f_c(x_k)$,   there are some  $x_k\to \bar{x}$ and $v_k\to v$ such that
\[
(v_k-\bar{v}_k, \pi_c(x_k)-x_k)\in \widehat{N}_{\gph F} (x_k, F(x_k))\quad \forall k,
\]
where  $\bar{v}_k:=F(x_k)- c (x_k-\pi_c(x_k))$.  In view of the fact that  $F$ and $\pi_c$  are locally Lipschitz continuous functions (cf. Lemma \ref{lem-basic} {\bf (f)}), we have $ \bar{v}_k  \to F(\bar{x})-c (\bar{x}-\pi_c(\bar{x}))$, $ x_k -\pi_c(x_k)\to  \bar{x} -\pi_c(\bar{x})$, and hence
$
(v-F(\bar{x})+c (\bar{x}-\pi_c(\bar{x})), \pi_c(\bar{x})-\bar{x})\in N_{\gph F} (\bar{x}, F(\bar{x}))$,  or in other words, $ v-F(\bar{x})+c (\bar{x}-\pi_c(\bar{x}))\in D^{*}F(\bar{x})(\bar{x}-\pi_c(\bar{x}))$.
This verifies that $v\in U$ and hence that  $\partial f_c(\bar{x})\subset U$.

To show $U \subset \partial f_c(\bar{x})$,  let   $v\in \displaystyle \left(D^{*}F(\bar{x})-c I\right)\left( \bar{x}-\pi_c(\bar{x}) \right)+F(\bar{x})$. Then we have
\[
 z:=v+c( \bar{x}-\pi_c(\bar{x}))- F(\bar{x})\in D^{*}F(\bar{x})(\bar{x}-\pi_c(\bar{x})) \Longleftrightarrow (z, -\bar{x}+\pi_c(\bar{x}))\in N_{\gph F}(\bar{x},F(\bar{x})).
\] 
 According to the definition of normal cone  (cf. Definition \ref{def-jihe-C}) and the definition of regular coderivative (cf. Definition \ref{def-set-valued-derivative}),   there exist $x_k\to \bar{x}$, $z_k\to z$ and $w_k\to \bar{x}-\pi_c(\bar{x})$ such that for all $k$,
\[
(z_k, -w_k)\in \widehat{N}_{\gph F}(x_k, F(x_k)) \Longleftrightarrow (z_k, -w_k)\in (\gph DF(x_k))^*,
\]
or explicitly,
\begin{equation}\label{nice-formu}
\langle z_k, w\rangle-\langle x_k-\pi_c(x_k), z\rangle\leq \langle w_k-x_k+\pi_c(x_k), z\rangle\quad \forall z \in DF(x_k)(w).
\end{equation}
By the Cauchy-Schwarz inequality and Lemma \ref{lem-lip-F-pro} {\bf (c)}, we have  for all $k$,
$$
\langle w_k-x_k+\pi_c(x_k), z\rangle\leq \epsilon_k \|w\|\quad \forall z \in DF(x_k)(w),
$$
where $\epsilon_k:={\rm lip} F(x_k)\| w_k-x_k+\pi_c(x_k)\|$. It then follows from (\ref{nice-formu}) that for all $k$,
$$
\langle z_k, w\rangle\leq \min \langle DF(x_k)(w), x_k-\pi_c(x_k)\rangle+\epsilon_k\|w\|\quad \forall w\in \R^n.
$$
By the formula for the subderivative $df_c(x_k)(w)$, we have for all $k$,
\begin{equation}\label{yjfh}
\langle z_k-c(x_k-\pi_c(x_k))+F(x_k),\; w\rangle\leq  df_c(x_k)(w)+\epsilon_k\|w\|\quad \forall w\in \R^n.
\end{equation}
In view of the fact that  $F$ and  $\pi_c$  are locally Lipschitz continuous functions (cf. Lemma \ref{lem-basic} {\bf (f)}) and by letting $k\to +\infty$, we have
$z_k-c(x_k-\pi_c(x_k))+F(x_k)\to z-c(\bar{x}-\pi_c(\bar{x}))+F(\bar{x})=v$, and $\epsilon_k\to 0$ (due to $\lip F(\cdot)$ being upper semicontinuous (\cite[Theorem 9.2]{roc}) and $w_k-x_k+\pi_c(x_k)\to 0$). Then by \cite[Proposition 10.46]{roc} and (\ref{yjfh}), we have $v\in \partial f_c(\bar{x})$. This completes the proof. \hfill $\Box$

By virtue of the formula for the limiting subdifferential set $\partial f_{c}(\bar{x})$ in Proposition \ref{subderivative000}, we can easily get the formula for the Clarke subdifferential set $\overline{\partial} f_{c}(\bar{x})$, which has been obtained first in \cite[Lemma 3.2]{xu}.

 \begin{corollary}\label{coro-gamma}
Let $\bar{x}\in \R^n$. We have
   $$
   \overline{\partial}  f_c(\bar{x})= \displaystyle \left(\overline{\partial} F(\bar{x})^T-cI\right) (\bar{x}-\pi_c(\bar{x})) +F(\bar{x}),
   $$
   where $\overline{\partial} F(\bar{x})$ denotes the generalized Jacobian   of $F$ at $\bar{x}$ (cf. (\ref{g-jac})).
  \end{corollary}
\noindent{\bf Proof.}  By Lemma \ref{lem-basic} {\bf (f)} and Lemma \ref{lem-lip-f-pro} {\bf (c)},   $f_c$ is locally Lipschitz continuous and hence
  $\overline{\partial} f_c(\bar{x})=\co (\partial f_c(\bar{x}))$.   The formula for $\overline{\partial} f_{c}(\bar{x})$ then  follows directly
   from  Proposition \ref{subderivative000}  and  the coderivative duality (\ref{code-dual}).  This completes the proof. \hfill $\Box$

\subsection{Subderivatives and subgradients of $f_{ab}$}
In parallel fashion as we have done in subsection \ref{sec-gamma}, we present in this subsection some differential properties of the D-gap function $f_{ab}$. Most of the proofs are omitted because they are very similar with the corresponding ones in subsection \ref{sec-gamma}.

\begin{proposition}\label{prop-sub-diff}
Let $\bar{x}\in \R^n$ and $w\in \R^n$. We have the following formulas:
\[
\begin{array}{ll}
&df_{ab}(\bar{x})(w)=\displaystyle (b-a)\langle \bar{x}-\pi_a(\bar{x}),\;w\rangle +\min \langle \left( DF(\bar{x})-bI \right)w, \; \pi_b(\bar{x})-\pi_a(\bar{x})\rangle,\\[0.1cm]
&\widehat{\partial}f_{ab}(\bar{x})=\left(\widehat{D}^{*}F(\bar{x})-bI\right)\left(\pi_b(\bar{x})-\pi_a(\bar{x})\right)+(b-a)(\bar{x}-\pi_a(\bar{x})),\\[0.1cm]
&\partial f_{ab}(\bar{x})=\left(D^{*}F(\bar{x})-bI\right)\left(\pi_b(\bar{x})-\pi_a(\bar{x})\right)+(b-a)(\bar{x}-\pi_a(\bar{x})),
\end{array}
\]
where 
$$
\min \langle \left( DF(\bar{x})-bI \right)w, \; \pi_b(\bar{x})-\pi_a(\bar{x})\rangle:=\min_{v\in DF(\bar{x})(w)}\langle \left( v-bw \right), \; \pi_b(\bar{x})-\pi_a(\bar{x})\rangle.
$$
\end{proposition}
\noindent{\bf Proof.} In view of the fact that $f_{ab}=f_a-f_b$ is a locally Lipschitz continuous function, we have
\[
df_{ab}(\bar{x})(w)=\liminf_{t\downarrow 0}\left[\frac{f_a(\bar{x}+t w)-f_a(\bar{x})}{t}-\frac{f_b(\bar{x}+t w)-f_b(\bar{x})}{t}\right].
\]
According to Lemma \ref{lem-basic} {\bf (g)}, we have for all $t$, $
f_a(\bar{x})\geq \langle F(\bar{x}),\bar{x}-\pi_a(\bar{x}+t w)\rangle-\frac{a}{2} ||\bar{x}-\pi_a(\bar{x}+tw)||^2$ and $
f_b(\bar{x}+t w)\geq \langle F(\bar{x}+tw),\bar{x}+tw-\pi_b(\bar{x})\rangle-\frac{b}{2} ||\bar{x}+tw-\pi_b(\bar{x})||^2$.
This, together with (\ref{chouxiang}) and the fact that $\pi_a$ and $\pi_b$ are locally Lipschitz continuous functions (see Lemma \ref{lem-basic} {\bf (f)}), implies that
\begin{eqnarray}
df_{ab}(\bar{x})(w)&\leq& \liminf_{t\downarrow 0}\langle \frac{F(\bar{x}+tw)-F(\bar{x})}{t}, \pi_b(\bar{x})-\pi_a(\bar{x}+tw)\rangle\nonumber\\[0.1cm]
&&-\lim_{t\downarrow 0}\frac{a}{2}\frac{||\bar{x}+tw-\pi_a(\bar{x}+tw)||^2-||\bar{x}-\pi_a(\bar{x}+tw)||^2}{t}\nonumber\\[0.1cm]
&&+\lim_{t\downarrow 0}\frac{b}{2}\frac{||\bar{x}+tw-\pi_b(\bar{x})||^2-||\bar{x}-\pi_b(\bar{x})||^2}{t}\nonumber\\[0.1cm]
&=&\min_{v\in DF(\bar{x})(w)}\langle v, \pi_b(\bar{x})-\pi_a(\bar{x})\rangle+\langle b (\bar{x}-\pi_b(\bar{x}))-a (\bar{x}-\pi_a(\bar{x})),w\rangle.\nonumber
\end{eqnarray}
To prove the inequality in the other direction, we simply follow a similar way by observing  from Lemma \ref{lem-basic} {\bf (g)} that for all $t$, $f_a(\bar{x}+t v)\geq  \langle F(\bar{x}+tv),\bar{x}+tv-\pi_a(\bar{x})\rangle-\frac{a}{2} ||\bar{x}+tv-\pi_a(\bar{x})||^2$  and $
f_b(\bar{x})\geq \langle F(\bar{x}),\bar{x}-\pi_b(\bar{x}+t v)\rangle-\frac{b}{2} ||\bar{x}-\pi_b(\bar{x}+tv)||^2$.
This completes the proof of the formula for $df_{ab}(\bar{x})(w)$. The other two formulas can be obtained in a similar way as we have done in Proposition \ref{subderivative000}.
\hfill $\Box$

 \begin{corollary}\label{coro-ab}
 Let $\bar{x}\in \R^n$. The following properties hold:
 \begin{description}
 \item[{\bf (a)}] We have the  formula for the Clarke subdifferential set of $f_{ab}$ at $\bar{x}$ as follows:
   $$
   \overline{\partial}  f_{ab}(\bar{x})= \displaystyle \left(\overline{\partial} F(\bar{x})^T-bI\right) (\pi_b(\bar{x})-\pi_a(\bar{x}))+(b-a)
(\bar{x}-\pi_a(\bar{x})).
   $$
     \item[{\bf (b)}] $\bar{x}$ solves  (VIP)  if and only if
 $ 0\in \partial f_{ab}(\bar{x})$
  and $\pi_a(\bar{x})=\pi_b(\bar{x})$.
   \end{description}
  \end{corollary}

  \begin{remark}
  The formula for $\overline{\partial} f_{ab}(\bar{x})$ was first  obtained in \cite[Lemma 3.3]{xu}, and then in \cite[Theorem 4.1]{ng2} and \cite[Theorem 3.1]{lig} for some generalized D-gap functions.
According to the generalized Fermat's rule \cite[Theorem 10.1]{roc}, the condition
  \begin{equation}\label{2ne}
   0\in \partial f_{ab}(\bar{x})
  \end{equation}
  is necessary for $\bar{x}$  to be  locally optimal for the  optimization  problem
$$
\min f_{ab}(x)\quad \mbox{s.t.}\quad x\in \R^n,
$$
and hence necessary for $\bar{x}$  to be a solution of (VIP) (cf.  Lemma \ref{lem-basic} {\bf (c)}). Another necessary condition for $\bar{x}$ to be a solution of (VIP) is, by  Lemma \ref{lem-basic} {\bf (c)},  the  equality
  \begin{equation}\label{3ne}
\pi_a(\bar{x})=\pi_b(\bar{x}).
  \end{equation}
Although these two necessary conditions  together become sufficient for $\bar{x}$  to be a solution of (VIP),
it is interesting to note that  either one alone  is not sufficient.

To see that (\ref{2ne}) alone is  not enough to guarantee that $\bar{x}$ solves (VIP), we simply consider the case that  $K=\R^n$ and $F$ is smooth with $\nabla F(\bar{x})^T F(\bar{x})=0$ but $F(\bar{x})\not=0$, for which case, (\ref{2ne}) holds as  $f_{ab}$ is  smooth with $\nabla f_{ab}(\bar{x})=\frac{b-a}{ab}\nabla F(\bar{x})^T F(\bar{x})=0$, but $\bar{x}$ does not solve (VIP) as $F(\bar{x})\not=0$. In this case,  (\ref{3ne}) does not hold as it amount to $F(\bar{x})=0$.

To see that (\ref{3ne}) alone  is not enough to guarantee that $\bar{x}$ solves (VIP), we simply consider the case that  $K=\R^n_+$ and $\bar{x}\in \R^n$ with $F_i(\bar{x})\geq 0$ and $\bar{x}_i<0$ for all $i$,   for which case,  (\ref{3ne}) holds as $\pi_a(\bar{x})=\pi_b(\bar{x})=0$,  but $\bar{x}$ does not solve (VIP) as $\bar{x}\not\in K$. In this case, (\ref{2ne}) does not hold as $0\not\in \partial f_{ab}(\bar{x})=\{(b-a)\bar{x}\}$.

It was shown in \cite[Theorem 4.3]{lig} that $\bar{x}$ solves (VIP) if and only if $0\in \overline{\partial} f_{ab}(\bar{x})$ and
\begin{equation}\label{fpling}
\left.
\begin{array}{c}
w\in T_{ab}(x, F, K),\quad Z\in \overline{\partial} F(x)\\[0.25cm]
Z^Tw\in T_{ab}(x, F, K)^*
\end{array}
\right\}\Rightarrow F(x)^Tw=0,
\end{equation}
where  $T_{ab}(x,F,K)$  is a cone defined as in Lemma \ref{lem-basic} {\bf (e)}. However, by resorting to Corollary \ref{coro-ab} {\bf (b)} and  noting that $\overline{\partial} f_{ab}(\bar{x})=\partial f_{ab}(\bar{x})$ in the presence of (\ref{3ne}),
we can refine \cite[Theorem 4.3]{lig}   as follows:   $\bar{x}$ solves (VIP) if and only if $0\in \overline{\partial} f_{ab}(\bar{x})$ and (\ref{3ne}) holds. Note that $\pi_a(\bar{x})$ and $\pi_b(\bar{x})$ are involved in the definition of $T_{ab}(x, F, K)$. So in contrast to the verification of (\ref{fpling}),  it is much easier to verify (\ref{3ne}). It is also  noteworthy that (\ref{3ne}) is implied by (\ref{2ne})  whenever the inequality
\begin{equation}\label{manman}
d(0,\partial f_{ab}(\bar{x}))\geq \mu\|\pi_b(\bar{x})-\pi_a(\bar{x})\|
\end{equation}
holds for some $\mu>0$. Inequalities in the form of (\ref{manman}) will play a crucial role in the next section.
\end{remark}

\section{ The Kurdyka-{\L}ojasiewicz  inequality and error bounds of $f_{ab}$} \label{sec-kl-error}
In this section, we study the KL inequality and error bounds for the D-gap function $f_{ab}$ by virtue of the formula for the limiting subdifferential sets   $\partial f_{ab}(x)$ presented in last section. Before summarizing our main results in Theorem \ref{theo-zuizhongyao}, we present in Lemmas \ref{lem-xys}-\ref{lem-suff} several results on necessary and sufficient conditions for the following inequalities:
$$
d(0,\partial f_{ab}(x))\geq \mu\|\pi_b(x)-\pi_a(x)\|\quad \forall x\in V,
$$
where $V$ is some open set in $\R^n$.

\begin{lemma}\label{lem-xys}
Let $x\in \R^n$ and let $\mu>0$. If $
 d(0, \partial f_{ab}(x))\geq \mu \|\pi_b(x)-\pi_a(x)\|$,
then
 \begin{equation}\label{chuanchuan}
 d(0, \partial f_{ab}(x) )\geq \frac{\mu(b-a)}{\mu+b+\lip F(x)}\|x-\pi_a(x)\|.
 \end{equation}
\end{lemma}
\noindent{\bf Proof.}  Let $w:=\pi_b(x)-\pi_a(x)$ and let $u:=x-\pi_a(x)$. By invoking the formula  for $\partial f_{ab}(x)$ in Proposition  \ref{prop-sub-diff},  we can find some $z^*\in D^*F(x)(w)$ such that $d(0,  \partial  f_{ab}(x))=\|z^*-bw+(b-a)u\|$. Then we get (\ref{chuanchuan}), as  we have
\[
\begin{array}{lll}
d(0,\partial f_{ab}(x))
&\geq& -\|z^*\|-b\|w\|+(b-a)\|u\| \\[0.1cm]
&\geq&-( b+\lip F(x) )\|w\|+(b-a)\|u\| \\[0.1cm]
&\geq&-\frac{b+\lip F(x)}{\mu} d(0,\partial f_{ab}(x))+(b-a)\|u\|,
\end{array}
\]
where the first   inequality follows from the triangle inequality, the second  one   from Lemma \ref{lem-lip-F-pro} {\bf (d)},
and the last one from the assumption that $ d(0, \partial f_{ab}(x))\geq \mu \|w\|$. This completes the proof.  \hfill $\Box$


\begin{lemma}\label{lem-KL}
Assume that $\lip F(x)$ is bounded from above on a nonempty subset $V$ of $\R^n$, as is true in particular when $V$ is  bounded. Then the  following properties are equivalent:
\begin{description}
      \item[{\bf (a)}] There is some $\mu>0$  such that
     $
d(0, \partial f_{ab}(x))\geq \mu  \sqrt{f_{ab}(x)}\quad \forall x\in V.
 $
  \item[{\bf (b)}]  There is some $\mu>0$   such that
 $
d(0, \partial f_{ab}(x))\geq \mu \|x-\pi_a(x)\| \quad \forall x\in V.
 $
  \item[{\bf (c)}]  There is some $\mu>0$  such that
 $
d(0, \partial f_{ab}(x))\geq \mu \|\pi_b(x)-\pi_a(x)\| \quad \forall x\in V.
 $
\end{description}
Therefore,    $f_{ab}$ satisfies the KL inequality  at any solution $\bar{x}$ of (VIP) with  an exponent of  $\frac{1}{2}$ if and only if  any of {\bf (a)}, {\bf (b)} and {\bf (c)} holds with $V$ being some neighborhood of $\bar{x}$.
\end{lemma}
\noindent{\bf Proof.}  The relations ${\bf (a)}   \Longleftrightarrow {\bf (b)}  \Longrightarrow {\bf (c)}$ follow directly from Lemma \ref{lem-basic} {\bf (a)}. As $\lip F(x)$ is upper semicontinuous (\cite[Theorem 9.2]{roc}), it follows from \cite[Corollary  1.10]{roc} that  $\lip F(x)$ is bounded from above  on each bounded subset of $\R^n$.  We now show  ${\bf (c)} \Longrightarrow {\bf (b)}$  by assuming that {\bf (c)} holds with some $\mu>0$ and that there is some $L>0$ such that
 $
\lip F(x)\leq L\;\forall x\in V.
$
By Lemma \ref{lem-xys}, we get  {\bf (b)} as we have
$$
d(0, \partial f_{ab}(x) )\geq \frac{\mu(b-a)}{\mu+b+\lip F(x)}\|x-\pi_a(x)\|\geq \frac{\mu(b-a)}{\mu+b+L}\|x-\pi_a(x)\|\quad\forall x\in V.
$$

Let $\bar{x}$ be a solution of (VIP). We first note that $f_{ab}$ is locally Lipschitz continuous with $f_{ab}\geq 0$ and $f_{ab}(\bar{x})=0$ (cf. Lemma \ref{lem-basic} {\bf (c)}). Then  $f_{ab}$ satisfies the KL inequality at $\bar{x}$ with an exponent of $\frac{1}{2}$ if, according to Definition \ref{def-kl},   {\bf (a)} holds with $V$ being some  bounded  neighborhood of $\bar{x}$.  By the previous argument,   {\bf (a)}, {\bf (b)} and {\bf (c)} are equivalent whenever $V$ is bounded, and therefore the last assertion is true. This completes the proof. \hfill $\Box$

  \begin{lemma}\label{lem-error-bound}
  Assume that the solution set of (VIP) is nonempty. If there are some $\mu\in (0, +\infty)$ and  $\varepsilon\in (0, +\infty]$ such that
   \begin{equation}\label{ruangang}
  \displaystyle d\left(0, \partial f_{ab}(x) \right)\geq \mu \|\pi_b(x)-\pi_a(x)\| \quad \forall x\in [f_{ab}< \varepsilon],
   \end{equation}
   and
\begin{equation}\label{ruangang000}
  L:=\displaystyle\sup_{x\in [0<f_{ab}< \varepsilon]}\lip F(x)<+\infty,
   \end{equation}
then
    \begin{equation}\label{stable-error-bound}
  \displaystyle\sqrt{\frac{b-a}{2}}\frac{\mu}{\mu+b+L}\,d\left(x,[f_{ab}\leq \theta]\right)\leq \left(\sqrt{f_{ab}(x)}-\sqrt{\theta}\right)_+\quad \forall \theta\in [0, \varepsilon),\; \forall x\in [f_{ab}< \varepsilon],
  \end{equation}
which,  in particular, implies the following error bound property:
    $$
  \displaystyle\sqrt{\frac{b-a}{2}}\frac{\mu}{\mu+b+L}\,d\left(x,[f_{ab}\leq 0]\right)\leq  \sqrt{f_{ab}(x)} \quad \forall  x\in [f_{ab}\leq \varepsilon].
  $$
  \end{lemma}
\noindent{\bf Proof.}  It suffices to show (\ref{stable-error-bound}) by assuming (\ref{ruangang}) and (\ref{ruangang000}) for some  given $\mu\in (0, +\infty)$ and  $\varepsilon\in (0, +\infty]$.  As the solution set of (VIP) is nonempty, we deduce from Lemma \ref{lem-basic} {\bf (c)} that
   $[f_{ab}\leq 0]\not=\emptyset$.  In what follows, we assume that $[0<f_{ab}<\varepsilon]$ is nonempty, for otherwise (\ref{stable-error-bound}) holds trivially.   Fix any $x\in [0<f_{ab}< \varepsilon]$.
  In view of (\ref{ruangang}) and (\ref{ruangang000}), we get from Lemma \ref{lem-xys} that
 $
    d(0, \partial f_{ab}(x) )\geq \displaystyle \frac{\mu(b-a)}{\mu+b+L}\|x-\pi_a(x)\|$.
 Then by Lemma \ref{lem-basic} {\bf (a)}, we  have
  $
    d(0, \partial f_{ab}(x) )\geq \displaystyle \frac{\mu\sqrt{2(b-a)}}{\mu+b+L}\sqrt{f_{ab}(x)}$.
 By some direct calculation,  we have $\displaystyle \partial \sqrt{f_{ab}}(x)=\frac{\partial f_{ab}(x)}{2\sqrt{f_{ab}}(x)}$
 and hence $d\left(0, \partial \sqrt{f_{ab}}(x)\right)\geq   \displaystyle\sqrt{\frac{b-a}{2}}\frac{\mu}{\mu+b+L}$.
Then by \cite[Lemma 2.1 {\bf (ii')}]{my12}, we have $$|\nabla \sqrt{f_{ab}} |(x)\geq  \displaystyle\sqrt{\frac{b-a}{2}}\frac{\mu}{\mu+b+L},$$
 where for a function $f:\R^n\to \R$ and a point $\bar{y}\in \R^n$,
 \[
|\nabla f|(\bar{y}):=\displaystyle\limsup_{y\rightarrow \bar{y},\;y\not=\bar{y}}\displaystyle\frac{(f(\bar{y})-f(y))_+}{\|y- \bar{y}\|}
\]
 denotes the the strong slope of $f$ at $\bar{y}$, introduced by De Giorgi et al.  \cite{gmt80}.  As $x\in [0<f_{ab}< \varepsilon]$ is chosen arbitrarily, we can apply \cite[Theorem 2.1]{ac04} to deduce that
 $$
 \begin{array}{ll}
 \displaystyle\inf_{0\leq \sqrt{\theta} <\sqrt{\varepsilon}}\inf_{x\in [\sqrt{\theta}<\sqrt{f_{ab}}<\sqrt{\varepsilon}\,]} \frac{\sqrt{f_{ab}(x)}-\sqrt{\theta}}{d\left(x, \left[\sqrt{f_{ab}}\leq \sqrt{\theta}\right] \right)}&= \displaystyle\inf_{x\in \left[0<\sqrt{f_{ab}}<\sqrt{\varepsilon} \right]} |\nabla \sqrt{f_{ab}} |(x)\\[0.5cm]
 & \displaystyle\geq \sqrt{\frac{b-a}{2}}\frac{\mu}{\mu+b+L},
 \end{array}
 $$
 from which, (\ref{stable-error-bound}) follows readily.   This completes the proof. \hfill $\Box$

Many existing conditions in the literature are sufficient for Lemma \ref{lem-KL} {\bf (c)} or  (\ref{ruangang}),  as can be seen from the following lemma, where we also provide  a new sufficient condition which can be considered as some restricted strong monotonicity.
\begin{lemma}\label{lem-suff}
Let $\mu>0$ and let $V\subset \R^n$ be open. Consider the following properties:
\begin{description}
  \item[{\bf (a)}]  $F$ is  strongly monotone on $V$ with modulus $\mu$, which holds in the case of $V$ being convex if and only if the following inequality holds  for all $x\in V$ where $F$ is differentiable:
          \begin{equation}\label{sm-mu}
     \displaystyle\langle \nabla F(x) w,\;w\rangle\geq \mu||w||^2\quad \forall w\in \R^n.
 \end{equation}
  \item[{\bf (b)}] The following holds for all $x\in V$ where $F$ is differentiable and $f_{ab}(x)>0$:
     $$
     \displaystyle\langle \nabla F(x) w,\;w\rangle\geq \mu||w||^2\quad \forall w\in T_{ab}(x,F,K).
 $$
\item[{\bf (c)}] The following holds for all $x\in V$ where $F$ is differentiable:
 $$\displaystyle\langle \nabla F(x) (\pi_a(x)-\pi_b(x)),\;\pi_a(x)-\pi_b(x)\rangle\geq \mu||\pi_a(x)-\pi_b(x)||^2.
 $$
\item[{\bf (d)}] $d(0, \partial f_{ab}(x) )\geq \mu \|\pi_b(x)-\pi_a(x)\| \quad \forall x\in V$.
\end{description}
We have
$
 {\bf (a)}   \Longrightarrow  {\bf (b)}  \Longrightarrow  {\bf (c)} \Longrightarrow  {\bf (d)}. 
$
\end{lemma}
\noindent{\bf Proof.}  According to \cite[Proposition 2.3 {\bf (b)}]{jiang},  the following holds for all $x\in V$:
       \begin{equation}\label{sm-mu000}
     \displaystyle\langle  Z w,\;w\rangle\geq \mu||w||^2\quad \forall Z\in \overline{\nabla}F(x),\; \forall w\in \R^n,
 \end{equation}
 if $F$ is  strongly monotone on $V$ with modulus $\mu$, and the converse is true whenever   $V$ is convex.
As $\nabla F(x)\in \overline{\nabla}F(x)$ when $F$ is differentiable at $x$,  (\ref{sm-mu}) is  implied by (\ref{sm-mu000}). To show   that (\ref{sm-mu000}) is implied by (\ref{sm-mu}), let $x\in V$ and let $Z\in \overline{\nabla}F(x)$. By the definition of $\overline{\nabla} F(x)$ (cf. (\ref{def-B-jac})), there is $x_k\to x$ such that $F$ is differentiable at $x_k$ for all $k$ and $\nabla F(x_k)\to Z$. Then by (\ref{sm-mu}),  we have for all sufficiently large $k$: $
\displaystyle\langle \nabla F(x_k) w,\;w\rangle\geq \mu||w||^2\quad \forall w\in \R^n$,
which implies (\ref{sm-mu000}) by letting $k \to \infty$.

By the previous argument,  we get {\bf (b)} from {\bf (a)} in a straightforward way. To get {\bf (c)} from {\bf (b)}, it suffices to note  the following facts: (1) $\pi_a(x)-\pi_b(x)\in T_{ab}(x,F,K)$ (cf. Lemma \ref{lem-basic} {\bf (e)}); (2) $\pi_a(x)=\pi_b(x)$ whenever $f_{ab}(x)=0$ (cf. Lemma \ref{lem-basic} {\bf (c)}).

We now  show ${\bf (c)} \Longrightarrow  {\bf (d)}$. Let $x\in V$. Set $w:=\pi_b(x)-\pi_a(x)$ and   $u:=x-\pi_a(x)$.  We first claim that the following holds for all $z^*\in \co D^*F(x)(w)$:
\begin{equation}\label{dataizi}
\langle  z^*, w \rangle\geqslant \mu ||w ||^2.
\end{equation}
By the coderivative duality (\ref{code-dual}) for a locally Lipschitz continuous mapping, we have $z^*\in \{A^Tw\mid A\in \co \overline{\nabla} F(x)\}$. Then there exist a positive integer $r$ and some $A^i\in \overline{\nabla} F(x)$    such that
\begin{equation}\label{convex hull}
z^*=\left(\sum_{i=1}^r \lambda^i A^i\right)^T w=\sum_{i=1}^r \lambda^i \left(A^i\right)^T w,
\end{equation}
where $\lambda^i\geq 0$ for all $i$ and $\sum_{i=1}^r \lambda^i=1$.
For each $A^i\in \overline{\nabla} F(x)$,  there exists by definition some sequence $\{x_k^i\}$  such that   $F$ is differentiable at $x_k^i$ for all $k$,
$x_k^i\to x$  and $\nabla F(x_k^i)\to A^i$ as $k\to \infty$. Then by {\bf (c)}, we have for all $k$ large enough,
\[
\langle  \nabla F(x_k^i) (\pi_a(x_k^i)-\pi_b(x_k^i)), \pi_a(x_k^i)-\pi_b(x_k^i) \rangle \geqslant  \mu ||\pi_b(x_k^i)-\pi_a(x_k^i) ||^2.
\]
Thus, by  noting that $\pi_a$ and $\pi_b$ are locally Lipschitz continuous and letting $k \to \infty$, we get
 $\langle  A^i (\pi_a(x)-\pi_b(x)), \pi_a(x)-\pi_b(x) \rangle\geqslant \mu ||\pi_b(x)-\pi_a(x) ||^2$, or in terms of $w$, $\langle (A^i)^Tw, w\rangle \geq \mu\|w\|^2$.
 This, together with (\ref{convex hull}), yields (\ref{dataizi}).

 By invoking the formula  for  $\partial f_{ab}(x)$ in Proposition \ref{prop-sub-diff},
we can find some  $\bar{z}^*\in D^*F(x)(w)\subset \co D^*F(x)(w)$ such that $d(0, \partial f_{ab}(x))=\|\bar{z}^*-bw+(b-a)u\|$.
Then we get {\bf (d)}, as we have
$
d(0, \partial f_{ab}(x))\,\|w\|\geq  \langle \bar{z}^*-bw+(b-a)u, \;w\rangle\geq \langle \bar{z}^*,   w\rangle\geq \mu\|w\|^2$,
 where the first inequality follows from the Cauchy-Schwarz inequality, the second one from  Lemma \ref{lem-basic} {\bf (d)}, and the last one from (\ref{dataizi}).
  This completes the proof. \hfill $\Box$

  \begin{remark}\label{rem-lig-ng}
   As $\nabla F(x)\in \overline{\nabla}F(x)\subset \overline{\partial} F(x)$ when $F$ is differentiable at $x$, Lemma  \ref{lem-suff} {\bf (b)} holds if the following holds for all $x\in V$ with $f_{ab}(x)>0$:
     \begin{equation}\label{ligy-ng}
  \langle Z^Tw, w\rangle\geq \mu ||w||^2\quad \forall Z\in \overline{\partial}F(x),\; \forall w\in T_{ab}(x,F,K).
  \end{equation}
 When  $V=\R^n$, the supremum of all  possible positive $\mu$ satisfying (\ref{ligy-ng}) can be reformulated  as
\begin{equation}\label{ligy-ng-mu}
\mu_{ab}:=\displaystyle\inf\{w^TZw \mid Z\in  \overline{\partial} F(x),\, w\in T_{ab}(x,F,K),\, \|w\|=1,\, f_{ab}(x)>0\}.
\end{equation}
The quantity $\mu_{ab}$ was first introduced for a general case  in \cite[Theorem 4.2]{lig}, where the condition $\mu_{ab}>0$ was utilized to study
the local error bounds for  $f_{ab}$.
  \end{remark}

\begin{remark}\label{rem-good-suff}
Lemma  \ref{lem-suff} {\bf (c)} can be reformulated as
\begin{equation}\label{oligei}
\langle z^*, \pi_b(x)-\pi_a(x)\rangle\geq \mu||\pi_a(x)-\pi_b(x)||^2\quad\forall x\in V,\; z^*\in \co D^*F(x)(\pi_b(x)-\pi_a(x)),
\end{equation}
or
\begin{equation}\label{oligei000}
\langle z, \pi_a(x)-\pi_b(x)\rangle\geq \mu||\pi_a(x)-\pi_b(x)||^2\quad\forall x\in V,\; z\in \co D_*F(x)(\pi_a(x)-\pi_b(x)),
\end{equation}
where  $D_*F(x)$ stands for the strict derivative mapping of $F$ at $x$ (cf. (\ref{def-strict-der})).
As $$\nabla F(x)^T(\pi_b(x)-\pi_a(x))\in\co D^*F(x)(\pi_b(x)-\pi_a(x))$$ and $$\nabla F(x)(\pi_a(x)-\pi_b(x))\in\co D_*F(x)(\pi_a(x)-\pi_b(x))$$  whenever $F$ is differentiable at $x$ (cf. (\ref{code-dual}) and (\ref{stri-dual})),  Lemma  \ref{lem-suff} {\bf (c)}  is clearly implied by (\ref{oligei}) or (\ref{oligei000}).
In the proof of ${\bf (c)} \Longrightarrow  {\bf (d)}$ in Lemma \ref{lem-suff}, we have already shown that (\ref{oligei}) is implied by  Lemma  \ref{lem-suff} {\bf (c)}.  By the coderivative duality (\ref{stri-dual}) for a locally Lipschitz continuous mapping,  we can show in a similar way that (\ref{oligei000}) is  also implied by Lemma  \ref{lem-suff} {\bf (c)}.
\end{remark}

\begin{example}
Let $A\in \R^{n\times n}$ and $q\in \R^n$ be such that $q+\rge A\not=\{0\}$, where $\rge A$ denotes the range space of $A$.   Consider a (VIP) instance with $K=\R^n$ and $F(x)=Ax+q$.    In this case, to find a solution of (VIP) is to find a solution to the linear equation $Ax+q=0$, which exists if and only if $q\in \rge A$.
Clearly, $F$ is continuously differentiable on $\R^n$ with $\nabla F(\cdot)=A$, implying that $f_{ab}$ is continuously differentiable on $\R^n$.
By some direct computation  we have
$$
\pi_b(x)-\pi_a(x)=\frac{b-a}{ab}(Ax+q),\quad f_{ab}(x)=\frac{b-a}{2ab}\|Ax+q\|^2,
$$
and
$$
 \nabla f_{ab}(x)=\frac{b-a}{ab}A^T(Ax+q), \quad T_{ab}(x, F, K)=\{w\mid \langle Ax+q, w\rangle\leq 0\}.
$$
Then in the case of $V:=\R^n$,  Lemma \ref{lem-suff} {\bf (a)-(d)} can be  reduced  respectively  to the following:
\begin{description}
  \item[{\bf (a)}] $A-\mu I$ is  positive-semidefinite on $\R^n$.
 \item[{\bf (b)}] $A-\mu I$ is  positive-semidefinite  on at least one closed-half space  containing the origin  and hence  on the whole space  $\R^n$.

  (Therefore, {\bf (a)} and {\bf (b)} coincide, both of which implies that $A$ is positive-definite on $\R^n$ and that the linear equation $Ax+q=0$ has a unique solution.)
\item[{\bf (c)}] $A-\mu I$ is  positive-semidefinite on the linear subspace $\R \{q\}+\rge A$, which entails positive-semidefiniteness of
 $A^TAA-\mu A^TA$  on $\R^n$ and is equivalent to it when $q\in \rge A$. (The latter property can be fulfilled for a symmetric matrix $A$ if and only if   $A$ is   positive-semidefinite and  $0<\mu<\lambda_i$ with $\lambda_i$ being any positive eigenvalue  of $A$.)
\item[{\bf (d)}]  $AA^T-\mu^2 I$ is  positive-semidefinite on the linear subspace $\R \{q\}+\rge A$, which entails positive-semidefiniteness of $(A^TA)^2-\mu^2A^TA$ on $\R^n$ and is equivalent to it when $q\in \rge A$.  (The latter property can be fulfilled as long as $0<\mu\leq \sqrt{\lambda_i}$ with $\lambda_i$ being any positive eigenvalue  of $A^TA$.)
\end{description}

Therefore, in the case of $q\in \rge A$ with $A$ being  symmetric  and positive-semidefinite (but not positive-definite),  Lemma \ref{lem-suff} {\bf (a)-(b)} cannot hold, but  Lemma \ref{lem-suff} {\bf (c)} can   as long as  $0<\mu<\lambda_i$ with $\lambda_i$ being any positive eigenvalue  of $A$.    This demonstrates that   Lemma \ref{lem-suff} {\bf (c)} can be strictly weaker than Lemma \ref{lem-suff} {\bf (a)-(b)}.  While in the case of $q\in \rge A$ with $A$ being symmetric but not positive-semidefinite,   Lemma \ref{lem-suff} {\bf (c)} cannot hold, but Lemma \ref{lem-suff} {\bf (d)} can  as long as $\mu$ is less than or equal to
the square root of the smallest positive eigenvalue of $A^TA$.   This demonstrates that  Lemma \ref{lem-suff} {\bf (d)} can  be strictly weaker than Lemma \ref{lem-suff} {\bf (c)}.

\end{example}

\begin{theorem}\label{theo-zuizhongyao}
Assume that any of  {\bf (a)-(d)} in Lemma \ref{lem-suff}  holds with some $\mu>0$ and  $V=\R^n$. Then the following properties hold:
  \begin{description}
    \item[{\bf (a)}] $f_{ab}$ is a KL  function  with an exponent of $\frac{1}{2}$.
    \item[{\bf (b)}] If $F$ is coercive on $\R^n$, then the solution set of (VIP) is nonempty and compact, and  $\sqrt{f_{ab}}$ has a local error bound on $\R^n$, i.e., the following holds for any given $\varepsilon>0$:
   $$
  \displaystyle\sqrt{\frac{b-a}{2}}\frac{\mu}{\mu+b+L}\,d\left(x,\,[f_{ab}\leq 0]\,\right)\leq  \sqrt{f_{ab}(x)} \quad \forall  x\in [f_{ab}\leq \varepsilon].
  $$
where $L$ is any number such that $L\geq \lip F(x)$ for all $x\in [0<f_{ab}<\varepsilon]$.
    \item[{\bf (c)}] If the solution set of (VIP) is nonempty and  $F$ is globally Lipschitz continuous with a constant $L>0$, then $\sqrt{f_{ab}}$ has a global error bound on $\R^n$, i.e., the following holds:
           $$
  \displaystyle\sqrt{\frac{b-a}{2}}\frac{\mu}{\mu+b+L}\,d\left(x,[f_{ab}\leq 0]\right)\leq  \sqrt{f_{ab}(x)} \quad \forall  x\in \R^n.
  $$
  \end{description}
 \end{theorem}
\noindent{\bf Proof.}  For each   $x$ that is a solution of (VIP), it follows from Lemma \ref{lem-KL} that  $f_{ab}$ is a KL  function at  $x$  with an  exponent of $\frac{1}{2}$.
 For each $x$  that is not a solution of (VIP), we claim that  $0\not\in \partial f_{ab}(x)$ and hence
   $f_{ab}$ is a KL  function  at $x$ with an  exponent of $0$, for otherwise the inclusion $0\in  \partial f_{ab}(x)$, together with the equality  $\pi_a(x)=\pi_b(x)$ as can be guaranteed by Lemma \ref{lem-suff} {\bf (d)}, would imply  that $x$ is a solution of (VIP) (cf. Corollary \ref{coro-ab} {\bf (b)}).   As a whole
 $f_{ab}$ is indeed a KL  function  with an exponent of $\frac{1}{2}$. This verifies {\bf (a)}.

 To show {\bf (b)}, fix any  $\varepsilon>0$ and let $\bar{L}:=\sup_{x\in [0<f_{ab}<\varepsilon]}\lip F(x)$. By the coerciveness of $F$ on $\R^n$ (hence on $K$), the solution set of (VIP) is nonempty and compact (cf. \cite[Proposition 2.2.7]{fac}), and the level set $[f_{ab}\leq \varepsilon]$ is bounded (cf. \cite[Lemma 4.1]{lig}).   As $\lip F(x)$ is upper semicontinuous (cf. \cite[Theorem 9.2]{roc}), it follows from \cite[Corollary  1.10]{roc} that  $\lip F(x)$ is bounded from above  on each bounded subset of $\R^n$. So we have $\bar{L}<+\infty$. Then by Lemma \ref{lem-error-bound}, we get  {\bf (b)} in a straightforward way.

 To show {\bf (c)}, we apply  Lemma \ref{lem-error-bound} again  by noting that
 $$ \sup_{x\in [0<f_{ab}<+\infty]}\lip F(x)\leq L.$$
 This completes the proof. \hfill $\Box$

\begin{remark}
In the presence of Lemma \ref{lem-suff} {\bf (a)} with some $\mu>0$ and $V=\R^n$ (i.e., $F$ is  strongly monotone on $\R^n$ with modulus $\mu$),  it was pointed out by \cite[Remark 2.1 (ii)]{lig} that $F$ is coercive on $\R^n$.  In this case,  Theorem \ref{theo-zuizhongyao} {\bf (b)} holds without explicitly assuming coerciveness.  While in the presence of Lemma \ref{lem-suff} {\bf (b)} with $V=\R^n$ and some $\mu>0$,  Theorem \ref{theo-zuizhongyao} {\bf (b)} can be deduced from \cite[Theorem 4.2]{lig}(cf. Remark \ref{rem-lig-ng}).  To the best of our knowledge,  all the results in Theorem \ref{theo-zuizhongyao},  except for the mentioned ones,  are new.
\end{remark}

\begin{example}[\cite{lig}, Example 4.4]\label{ex-li-ng}  Consider a (VIP) instance with $K=\R_+^2$ and  $F: \R^2\to \R^2$ being given  by
$
F(x)=\left(x_1+(x_1)_+(x_2)_+, \quad x_2+\frac{3}{2}(x_1)_+\right)^T$.
Clearly,  $F$ is differentiable at $x\in \R^2$ if and only if $x_1x_2\not=0$, and moreover,
\begin{center}
$\nabla F(x)=\left\{
     \begin{array}{ll}
     \left(
       \begin{array}{cc}
         1+x_2 & x_1 \\
         \frac{3}{2} & 1 \\
       \end{array}
     \right)
      & \mbox{if}\; x_1> 0, x_2>0,\\[0.1cm]
     \left(
       \begin{array}{cc}
         1 & 0 \\
         \frac{3}{2} & 1 \\
       \end{array}
     \right)
      & \mbox{if}\; x_1> 0, x_2<0,\\[0.1cm]
      \left(
       \begin{array}{cc}
         1 & 0 \\
         0 & 1 \\
       \end{array}
     \right)
      & \mbox{if}\; x_1< 0, x_2\neq 0.
     \end{array}
     \right.
$\end{center}
Let $a\in(0,1)$ and $b=1$.  According to \cite[Example 4.4]{lig}, $F$ is coercive  and   not monotone on $\R^2$, and $\sqrt{f_{ab}}$ has a local error bound on $\R^2$ (with some error bound modulus expressed in an abstract way),  and $\mu_{ab}\geq 1$, where $\mu_{ab}$ is defined by (\ref{ligy-ng-mu}).

In what follows, by virtue of Lemma \ref{lem-suff} {\bf (c)},  we can show that $\mu_{ab}=1$ and that some error bound modulus expressed in an explicit way
can be provided. First, by some direct calculation,  we have
$\pi_b(x)=(0, 0)^T$ for all $x\in \R^2$
and
\begin{center}
$
\pi_a(x)-\pi_b(x)=\left\{
\begin{array}{ll}
\displaystyle\left(\frac{a-1}{a}x_1,\; 0\right)^T &\mbox{if}\;x_1\leq 0, x_2\geq 0,\\[0.1cm]
\displaystyle\left(\frac{a-1}{a}x_1,\; \frac{a-1}{a}x_2\right)^T &\mbox{if}\;x_1\leq 0, x_2\leq 0,\\[0.1cm]
\displaystyle\left(0,\; \frac{a-1}{a}x_2-\frac{3}{2a}x_1\right)^T &\mbox{if}\;0\leq x_1\leq \frac{2(a-1)}{3}x_2,\\[0.1cm]
\displaystyle\left(0,\; 0\right)^T &\mbox{otherwise}.\\
\end{array}
\right.
$\end{center}
Then it is straightforward to verify that the inequality
$$
\displaystyle\langle \nabla F(x) (\pi_a(x)-\pi_b(x)),\;\pi_a(x)-\pi_b(x)\rangle\geq \mu||\pi_a(x)-\pi_b(x)||^2
$$
holds for all $x\in \R^2$ with $x_1x_2\not=0$ if and only if  $0<\mu\leq 1$. That is,  Lemma \ref{lem-suff} {\bf (c)} holds with $V=\R^2$ if and only if $0<\mu\leq 1$.  As Lemma \ref{lem-suff} {\bf (c)} is implied by Lemma \ref{lem-suff} {\bf (b)},  we deduce that Lemma \ref{lem-suff} {\bf (b)} cannot hold with $V=\R^2$ and $\mu>1$,  which implies that $\mu_{ab}$ cannot be greater than 1 (cf. Remark \ref{rem-lig-ng}). Therefore, we confirm that  $\mu_{ab}=1$. Furthermore, we can apply Theorem \ref{theo-zuizhongyao} to get the following: (i) $f_{ab}$ is a KL  function  with an exponent of $\frac{1}{2}$; (ii) $\sqrt{f_{ab}}$ has a local error bound on $\R^2$, i.e.,  for any given $\varepsilon>0$,
   $$
  \displaystyle\sqrt{\frac{b-a}{2}}\frac{1}{1+b+L}\,d\left(x,[f_{ab}\leq 0]\right)\leq  \sqrt{f_{ab}(x)} \quad \forall  x\in [f_{ab}\leq \varepsilon],
  $$
where $L$ is any number such that $L\geq \sup_{x\in [0<f_{ab}<\varepsilon]}\lip F(x)$.
\end{example}

\section{A derivative free descent method for (VIP)}\label{sec-method}
In this section, we analyze the convergence behavior of the following descent algorithm with an Armijo line search, which is essentially the same as those studied in \cite{ka,lig,peng1,qu,xu,yam1}, especially the same in the way how descent directions are chosen.

\noindent{\bf Algorithm}
\begin{description}
  \item[Step 1.]   Set  $0<a<b$ and  $0<\rho<1$. Choose three positive constants $\alpha, \beta, \tau$ such that $\beta$ and $\tau$ are small   and that   $\alpha$ is close to $b-a$. Select a start point $x_0\in \R^n$, and set $n=0$.
       \item[Step 2.] If $f_{ab}(x_n)=0$,  stop. Otherwise,  go to Step 3.
        \item[Step 3.]  Let $u_n=\pi_{a}(x_n)-x_n$ and $w_n=\pi_{a}(x_n)-\pi_{b}(x_n)$. If $\beta ||u_n||< ||w_n||$,  set $d_n=w_n$ and select  $m_n$ as  the smallest nonnegative integer $m$ such that
\begin{equation}\label{iteration}
f_{ab}(x_n+\rho^m d_n)-f_{ab}(x_n)\leq -\tau \rho^m ||d_n||^2.
\end{equation}
Otherwise, set $d_n=u_n$ and  select  $m_n$ as the smallest nonnegative integer $m$ such that
\begin{equation}\label{iteration1}
f_{ab}(x_n+\rho^m d_n)-f_{ab}(x_n)\leq -\left(b-a-\alpha\right) \rho^m ||d_n||^2.
\end{equation}
         \item[Step 4.] Set $t_n=\rho^{m_n}$, $x_{n+1} = x_n + t_n d_n$ and $n = n+ 1$, and go to Step 2.
\end{description}

In what follows, we  make the following assumptions.
\begin{description}
\item[{\bf Assumption (i)}] The level set $[f_{ab}\leq f_{ab}(x_0)]$ is bounded, which can be guaranteed by the coerciveness of $F$ on $\R^n$ as pointed out by   \cite[Lemma 4.1]{lig}.
\item[{\bf Assumption (ii)}] $F$ is globally Lipschitz continuous   with a constant $L>0$ (implying that $f_{ab}$,  $\pi_a$ and $\pi_b$ are all globally Lipschitz continuous).
\item[{\bf Assumption (iii)}] There exists some $\mu^*>0$ such that  the inequality
$$
\displaystyle\langle \nabla F(x) (\pi_a(x)-\pi_b(x)),\;\pi_a(x)-\pi_b(x)\rangle\geq \mu^*||\pi_a(x)-\pi_b(x)||^2
$$
holds for all $x\in\R^n$ where  $F$ is differentiable. This implies by Theorem \ref{theo-zuizhongyao} that $f$ is a KL function with an exponent of $\frac12$, and by Remark \ref{rem-good-suff} and (\ref{txds-strict})  that
  \[
  \min_{z\in DF(x)(\pi_a(x)-\pi_b(x))}\langle z, \pi_a(x)-\pi_b(x)\rangle\geq \mu^*||\pi_a(x)-\pi_b(x)||^2\quad \forall x\in \R^n.
  \]
  \item[{\bf Assumption (iv)}] The parameters $\alpha, \beta, \tau$ in the Algorithm  are chosen such that
     \[
     0<\beta<\frac{b-a}{b+L},\quad  (b+L)\beta<\alpha< b-a,\quad 0<\tau<\mu^*.
     \]
\end{description}

To begin with, we give two technical lemmas, which are helpful for our further analysis.

 \begin{lemma}\label{lem-upper-bound}
 Under {\bf Assumption  (ii)}, we have
\[
\|v\|\leq (b+L)\|\pi_b(x)-\pi_a(x)\|+(b-a) \|x-\pi_a(x)\|\quad \forall x\in \R^n,\;\forall v\in \partial f_{ab}(x).
\]
\end{lemma}
\noindent{\bf Proof.} In view of Lemma \ref{lem-lip-F-pro} {\bf (d)}  and {\bf Assumption  (ii)}, we get this result directly from
  the formula for $\partial f_{ab}(x)$ presented in Proposition \ref{prop-sub-diff}.  The proof is completed. \hfill $\Box$

\begin{lemma}\label{mean-value theorem} Consider  a locally Lipschitz continuous function $g:\R^n\to \R$. For some $x\in \R^n$ and $w\in \R^n\backslash \{0\}$, assume that there are some $\sigma>0$ and  $0<t_0<t_1$ such that
$$g(x+t_0 w)-g(x)\leq -\sigma t_0 ||w||^2 \;\;\mbox{and}\;\; g(x+t_1 w)-g(x)> -\sigma t_1 ||w||^2.$$
Then there exist some $\theta^*\in (0, 1)$ and $v^*\in \partial g(x+\theta^*t_1 w)$ such that
$$g(x+t_1 w)-g(x)=t_1\langle v^*, w\rangle.$$
\end{lemma}
 \noindent{\bf Proof.} Define $\varphi:\R\rightarrow \R$ by $\varphi(\theta):=g(x+\theta t_1 w)-g(x)+\theta[g(x)-g(x+t_1w)]$. Clearly,  $\varphi$ is locally Lipschitz continuous, and  $\varphi(0)=\varphi(1)=0$.  Moreover, it follows from the assumption that $\varphi(t_0/t_1)=g(x+ t_0 w)-g(x)+ (t_0/t_1)[g(x)- g(x+t_1 w)]<0$.
This entails the existence of at least one $\theta^*\in (0,1)$ such that $\varphi$ attains its minimum over $[0,1]$ at $\theta^*$, implying  by the Fermat's rule  that  $0\in \partial \varphi(\theta^*)$.  In view of the local Lipschitzian continuity of $g$, we get from the calculus rules  \cite[Exercise 8.8 and Theorem 10.6]{roc} that $
\partial \varphi(\theta^*)\subset g(x)- g(x+t_1 w)+ \left\{t_1 \langle v, w\rangle\mid v\in \partial g(x+\theta^* t_1 w)\right\}$.
This completes the proof. \hfill $\Box$

\begin{proposition}\label{iteration2}
Under {\bf Assumptions  (ii)-(iv)},  Step 3 of  the Algorithm is well defined.
\end{proposition}
\noindent{\bf Proof.}  To show that Step 3 in the Algorithm is well defined, it suffices to show that  if $\beta ||u_n||< ||w_n||$,
$-d(-f_{ab})(x_n)(w_n)<-\tau \|w_n\|^2$, and if $\beta ||u_n||\geq ||w_n||$,   $-d(-f_{ab})(x_n)(u_n)<-(b-a-\alpha)\|u_n\|^2$.  Following from the proof of
the formula for   $df_{ab}(\bar{x})(w)$ in Proposition \ref{prop-sub-diff}, we get  the formula for the subderivative of $-f_{ab}$ at a point $\bar{x}\in \R^n$
   as follows:
$$
-d(-f_{ab})(\bar{x})(w)=\displaystyle (b-a)\langle \bar{x}-\pi_a(\bar{x}),\;w\rangle -\min \langle \left( DF(\bar{x})-bI \right)w, \; -\pi_b(\bar{x})+\pi_a(\bar{x})\rangle.
$$
In the case of $\beta ||u_n||< ||w_n||$, we have
\[
 \begin{array}{lll}
&&-d(-f_{ab})(x_n)(w_n) \\[0.1cm]
&=& \langle b (x_n-\pi_b(x_n))-a (x_n-\pi_a(x_n)),w_n\rangle-\min_{z\in DF(x_n)(w_n)}\langle z, w_n\rangle \\[0.1cm]
  &\leq& -\min_{z\in DF(x_n)(w_n)}\langle z, w_n\rangle\\[0.1cm]
  &\leq& -\mu^*||w_n||^2\\[0.1cm]
  &<&-\tau ||w_n||^2,
\end{array}
\]
where the first inequality follows from Lemma \ref{lem-basic} {\bf (d)}, the second inequality follows from  {\bf Assumption  (iii)}, and the third inequality follows from {\bf Assumption  (iv)}.  In the case of $\beta ||u_n||\geq ||w_n||$, we  have
\[
 \begin{array}{lll}
&&-d(-f_{ab})(x_n)(u_n) \\[0.1cm]
&=& \langle b (x_n-\pi_b(x_n))-a (x_n-\pi_a(x_n)),u_n\rangle-\min_{z\in DF(x_n)(u_n)}\langle z, w_n\rangle\\[0.1cm]
  &=& -(b-a)||u_n||^2+b \langle \pi_a(x_n)-\pi_b(x_n),u_n\rangle+\max_{z\in DF(x_n)(u_n)}\langle z, -w_n\rangle\\[0.1cm]
  &\leq& -[(b-a)-b\beta]||u_n||^2+\max_{z\in DF(x_n)(u_n)}\langle z, -w_n\rangle\\[0.1cm]
  &\leq&  -[(b-a)-b\beta]||u_n||^2+L||u_n||\cdot||w_n||\\[0.1cm]
  &\leq &-[(b-a)-(b+L)\beta]||u_n||^2\\[0.1cm]
  &<& -[(b-a)-\alpha]||u_n||^2,
\end{array}
\]
where the first  inequality follows by using  the Cauchy-Schwarz inequality and the inequality  $\beta ||u_n||\geq ||w_n||$, the second inequality follows from Lemma \ref{lem-lip-F-pro} {\bf (c)} and  {\bf Assumption (ii)},
  the third inequality follows from the inequality  $\beta ||u_n||\geq ||w_n||$,  and the last inequality follows from {\bf Assumption  (iv)}.  This completes the proof. \hfill $\Box$

\begin{proposition}\label{positive bound}
Assume that  the  sequence  $\{x_n\}$  generated by the Algorithm satisfies  $f_{ab}(x_n)>0$ for all $n$.
Under {\bf Assumptions  (ii)-(iv)},  there is some $t^*>0$ such that  $t_n\geq t^*$ for all $n$, i.e., the step length sequence  $\{t_n\}$ generated by the Algorithm has a lower bound.
\end{proposition}
\noindent{\bf Proof.}   Recall that in Step 3 of the Algorithm, we set $u_n:=\pi_a(x_n)-x_n$, $w_n:=\pi_a(x_n)-\pi_b(x_n)$, and $d_n:=u_n$ if $\beta\|u_n\|\geq \|w_n\|$, and $d_n:=w_n$ if $\beta\|u_n\|< \|w_n\|$.  In view of the setting for $d_n$ and our assumption that $f_{ab}(x_n)>0$ for all $n$, we get from  Lemma \ref{lem-basic} {\bf (c)}  that $d_n\not=0$ for all $n$.

Suppose by contradiction that the step length sequence  $\{t_n\}$ does not have a positive lower bound, i.e.,  by taking a subsequence if necessary we  assume that $t_n\to 0+$ as $n\to +\infty$.   Due to $t_n=\rho^{m_n}$, we have  $m_n\to +\infty$ as $n\to +\infty$. Without loss of generality, we may assume that $m_n\geq 1$ for all $n$.   In view of the line search strategy in Step 3 of the Algorithm, we apply Lemma \ref{mean-value theorem} to get
\begin{equation}\label{mean-value}
f_{ab}(x_n+\rho^{m_n-1}d_n)-f_{ab}(x_n)=\rho^{m_n-1}\langle v_n, d_n\rangle\quad \forall n,
\end{equation}
where  $v_n\in \partial f_{ab}(y_n)$ with $y_n:=x_n+\theta_n^*\rho^{m_n-1}d_n$ and $\theta_n^*\in (0, 1)$. By the formula for $\partial f_{ab}(y_n)$  in Proposition \ref{prop-sub-diff},  there exists some $z_n^*\in D^*F(\pi_b(y_{n})-\pi_a(y_{n}))$ such that
\begin{equation}\label{equality 0}
  v_n=z_n^*+b (y_{n}-\pi_b(y_{n}))-a (y_{n}-\pi_a(y_{n})).
\end{equation}
In view of   Lemma \ref{lem-lip-F-pro} {\bf (d)} and {\bf Assumption  (ii)}, we have
\begin{equation}\label{pudilan}
||z_{n}^*||\leq L ||\pi_b(y_{n})-\pi_a(y_{n})||.
\end{equation}

First, we consider the case that $\beta \|u_n\|\geq \|w_n\|$ in Step 3. In this case, we have  $d_n=u_n=\pi_a(x_n)-x_n$ and $y_n:=x_n+\theta_n^*\rho^{m_n-1}u_n$. Due  to the line search strategy proposed in the Algorithm, we have
$
f_{ab}(x_n+\rho^{m_n-1}u_n)-f_{ab}(x_n)>-(b-a-\alpha)\rho^{m_n-1}\|u_n\|^2$.
 This, together with (\ref{mean-value}), (\ref{equality 0})  and (\ref{pudilan}),  implies that
\begin{eqnarray}\label{inequality 3}
-(b-a-\alpha)||u_n||^2&&< \langle v_n, u_n\rangle \nonumber\\[0.1cm]
&&=\langle z_{n}^*, u_n\rangle+ \langle b (y_{n}-\pi_b(y_{n}))-a (y_{n}-\pi_a(y_{n})), u_n\rangle \nonumber\\[0.1cm]
&&= \langle z_{n}^*, u_n\rangle+ b\langle \pi_a(y_{n})-\pi_b(y_{n}), u_n\rangle+(b-a)\langle y_{n}-\pi_a(y_{n}), u_n\rangle\nonumber\\[0.1cm]
&&\leq  (L+b) ||\pi_b(y_{n})-\pi_a(y_{n})||\cdot||u_n||-(b-a)\langle \pi_a(y_{n})-y_n, u_n\rangle.
\end{eqnarray}
 Moreover, by  {\bf Assumption (ii)}, we have
\begin{eqnarray}\label{inequality 4.2}
||\pi_a(y_{n})-\pi_b(y_{n})||&\leq &||w_n||+||\pi_a(y_{n})-\pi_b(y_{n})-w_n||\nonumber \\[0.1cm]
&\leq &\|w_n\|+||\pi_a(y_{n})-\pi_a(x_n)\|+\|\pi_b(y_{n})-\pi_b(x_n)|| \nonumber\\[0.1cm]
&\leq&\beta\|u_n\|+(1+\frac{L}{a})\|y_n-x_n\|+(1+\frac{L}{b})\|y_n-x_n\| \nonumber\\[0.1cm]
&=& [\beta+(2+\frac{L}{a}+\frac{L}{b})\theta_n^*\rho^{m_n-1}]||u_n||,
\end{eqnarray}
and
\[
\begin{array}{lll}
||\pi_a(y_{n})-y_{n}-u_n||&= &||\pi_a(y_{n})-y_{n}-\pi_a(x_n)+x_n|| \\[0.1cm]
&\leq & \|\pi_a(y_n)-\pi_a(x_n)\|+\|y_n-x_n\|  \\[0.1cm]
&\leq&(2+\frac{L}{a})||y_n-x_n||=(2+\frac{L}{a})\theta_n^*\rho^{m_n-1}||u_n||.
\end{array}
\]
 The latter condition entails  that
\begin{equation}\label{ineqq}
\langle \pi_a(y_{n})-y_{n}, u_n\rangle=||u_n||^2+(2+\frac{L}{a})\theta_n^*\rho^{m_n-1}||u_n||^2\langle c_n, \frac{u_n}{\|u_n\|}\rangle,
\end{equation}
where $c_n:=\frac{\pi_a(y_{n})-y_{n}-u_n}{(2+\frac{L}{a})\theta_n^*\rho^{m_n-1}||u_n||}$
 having the property that  $\|c_n\|\leq 1$. Combining (\ref{inequality 3}-\ref{ineqq}), we have
\begin{equation}\label{key-ineq-1}
\begin{array}{lll}
-(b-a-\alpha)&<& (L+b) [\beta+(2+\frac{L}{a}+\frac{L}{b})\theta_n^*\rho^{m_n-1}] \\[0.1cm]
&&-(b-a)[1+(2+\frac{L}{a})\theta_n^*\rho^{m_n-1}\langle c_n, \frac{u_n}{||u_n||}\rangle].
\end{array}
\end{equation}

Next, we consider the case that $\beta \|u_n\|<\|w_n\|$ in Step 3. In this case, we have  $d_n=w_n=\pi_a(x_n)-\pi_b(x_n)$ and $y_n:=x_n+\theta_n^*\rho^{m_n-1}w_n$.   Due  to the line search strategy proposed in the Algorithm, we have $
f_{ab}(x_n+\rho^{m_n-1}w_n)-f_{ab}(x_n)>-\tau\rho^{m_n-1}\|w_n\|^2$, which, together with (\ref{mean-value}), (\ref{equality 0}) and  (\ref{pudilan}), implies that
\begin{equation}\label{yibudaowei}
\begin{array}{lll}
&&-\tau\|w_n\|^2\\[0.1cm]
&<& \langle v_n, \;w_n\rangle\\[0.1cm]
&=&\langle z_n^*+b(y_n-\pi_b(y_n))-a(y_n-\pi_a(y_n)), w_n\rangle\\[0.1cm]
&\leq & \langle z_n^*,\,\pi_a(y_n)-\pi_b(y_n)\rangle+ \langle z_{n}^*, w_n-(\pi_a(y_{n})-\pi_b(y_{n}))\rangle\\[0.1cm]
&&+ \langle b (y_{n}-\pi_b(y_{n}))-a (y_{n}-\pi_a(y_{n})), w_n-(\pi_a(y_{n})-\pi_b(y_{n}))\rangle\\[0.1cm]
&\leq & -\mu^* ||\pi_a(y_{n})-\pi_b(y_{n})||^2+ \langle z_{n}^*, w_n-(\pi_a(y_{n})-\pi_b(y_{n}))\rangle\\[0.1cm]
&&+ \langle b (y_{n}-\pi_b(y_{n}))-a (y_{n}-\pi_a(y_{n})), w_n-(\pi_a(y_{n})-\pi_b(y_{n}))\rangle\\[0.1cm]
&\leq& -\mu^* ||\pi_a(y_{n})-\pi_b(y_{n})||^2+L ||\pi_a(y_{n})-\pi_b(y_{n})||\cdot||w_n-(\pi_a(y_{n})-\pi_b(y_{n}))||\\[0.1cm]
&&+[(b-a)||\pi_a(y_{n})-y_{n}||+b||\pi_a(y_{n})-\pi_b(y_{n})||]||w_n-(\pi_a(y_{n})-\pi_b(y_{n}))||,
\end{array}
\end{equation}
where the second inequality follows from Lemma \ref{lem-basic} {\bf (d)}, the third one from {\bf Assumption  (iii)}, the last one from Cauchy-Schwarz inequality.
 Moreover, by  {\bf Assumption  (ii)}, we have
 \begin{equation}\label{inequality 2.2}
||\pi_a(y_{n})-\pi_b(y_{n})-w_n||\leq (2+\frac{L}{a}+\frac{L}{b})||y_n-x_n||=(2+\frac{L}{a}+\frac{L}{b})\theta_n^*\rho^{m_n-1}||w_n||,
\end{equation}
\begin{equation}\label{inequality 2.4}
||\pi_a(y_{n})-\pi_b(y_{n})||\leq [1+(2+\frac{L}{a}+\frac{L}{b})\theta_n^*\rho^{m_n-1}]||w_n||,
\end{equation}
\begin{eqnarray}\label{inequality 222}
 ||\pi_a(y_{n})-y_{n}||&\leq& \|u_n\|+||\pi_a(y_{n})-y_{n}-u_n||\nonumber\\[0.1cm]
 &\leq& ||u_n||+(2+\frac{L}{a})\theta_n^*\rho^{m_n-1}||w_n||\nonumber\\[0.1cm]
 &\leq& [\frac{1}{\beta}+(2+\frac{L}{a})\theta_n^*\rho^{m_n-1}]||w_n||
\end{eqnarray}
and then there exists $b_n$ with $||b_n||\leq 1$ such that
\begin{equation}\label{inequality 2.21}
\pi_a(y_{n})-\pi_b(y_{n})=w_n+(2+\frac{L}{a}+\frac{L}{b})\theta_n^*\rho^{m_n-1}||w_n||b_n.
\end{equation}
Combining (\ref{yibudaowei}-\ref{inequality 2.21}), we have
\begin{eqnarray}\label{key-ineq-2}
-\tau &<& -\mu^* [1+2\langle \frac{w_n}{||w_n||}, (2+\frac{L}{a}+\frac{L}{b})\theta_n^*\rho^{m_n-1} b_n\rangle+ (2+\frac{L}{a}+\frac{L}{b})^2(\theta_n^*\rho^{m_n-1})^2||b_n||^2] \nonumber\\[0.1cm]
&&+ L[1+(2+\frac{L}{a}+\frac{L}{b})\theta_n^*\rho^{m_n-1}](2+\frac{L}{a}+\frac{L}{b})\theta_n^*\rho^{m_n-1}\nonumber\\[0.1cm]
&&+ (b-a)[\frac{1}{\beta}+(2+\frac{L}{a})\theta_n^*\rho^{m_n-1}](2+\frac{L}{a}+\frac{L}{b})\theta_n^*\rho^{m_n-1}\nonumber\\[0.1cm]
&&+ b[1+(2+\frac{L}{a}+\frac{L}{b})\theta_n^*\rho^{m_n-1}](2+\frac{L}{a}+\frac{L}{b})\theta_n^*\rho^{m_n-1}.
\end{eqnarray}

Our assumption that $f_{ab}(x_n)>0$ for all $n$ suggests that there are infinitely many positive integers $n$ such that either $\beta \|u_n\|\geq \|w_n\|$ or
$\beta \|u_n\|<\|w_n\|$, implying that there are infinitely many positive integers $n$ such that  either the inequality (\ref{key-ineq-1}) or (\ref{key-ineq-2}) holds. In view of $\rho^{m_n-1}\rightarrow 0+$, we have  correspondingly either
$-(b-a-\alpha)\leq (L+b)\beta-(b-a)$ or $-\tau \leq  -\mu^*$, both contradicting to {\bf  Assumption (iv)}.  This contradiction  indicates that the step length sequence  $\{t_n\}$ generated by the Algorithm has a positive lower bound. This completes the proof. \hfill $\Box$

\begin{proposition}\label{prop-key conditions}
Assume that  the  sequence  $\{x_n\}$  generated by the Algorithm satisfies  $f_{ab}(x_n)>0$ for all $n$. Under  {\bf Assumptions  (ii)-(iv)},  the following inequalities hold for all $n$:
\begin{equation}\label{sufficient decrease condition}
  f_{ab}(x_{n+1})-f_{ab}(x_n)\leq -M_1||x_{n+1}-x_n||^2
\end{equation}
and
\begin{equation}\label{relative error condition}
  d(0, \partial f_{ab}(x_n))\leq  \frac{M_2}{t^*} ||x_{n+1}-x_n||,
\end{equation}
where $M_1:= \min\{b-a-\alpha, \tau\} $,  $M_2:=L+b+\frac{b-a}{\beta}$ and $t^*$ is a positive lower bound of   $\{t_n\}$.
\end{proposition}
\noindent{\bf Proof.}
By Steps 3 and 4 of the Algorithm, we have   $0<t_n\leq 1$,  $x_{n+1}=x_n+t_nd_n$ and
$f_{ab}(x_{n+1})-f_{ab}(x_n)\leq -M_1t_n\|d_n\|^2$ for all $n$, from which we get (\ref{sufficient decrease condition}) immediately.
 By Lemma \ref{lem-upper-bound}, we have
$
d(0, \partial f_{ab}(x_n))\leq (L+b)||w_n||+(b-a) ||u_n||$,
where $L$ is given as in {\bf Assumption  (ii)}, and  $w_n=\pi_a(x_n)-\pi_b(x_n)$ and $u_n=\pi_a(x_n)-x_n$ are set as in Step 3.
If $\beta ||u_n||<||w_n||$,  we get from  Steps 3 and  4 of the Algorithm that   $\|x_{n+1}-x_n\|=t_n\|w_n\|$ and hence that
\[
(L+b)||w_n||+(b-a) ||u_n||< (L+b+\frac{b-a}{\beta})||w_n||=\frac{M_2}{t_n} ||x_{n+1}-x_n||.
\]
Alternatively  if $\beta ||u_n||\geq ||w_n||$, we get from  Steps 3 and  4 of the Algorithm that $\|x_{n+1}-x_n\|=t_n\|u_n\|$ and hence that
\[
(L+b)||w_n||+(b-a) ||u_n||\leq \beta (L+b+\frac{b-a}{\beta}) ||u_n|| \leq \frac{M_2}{t_n} ||x_{n+1}-x_n||,
\]
where the second inequality follows from the fact that $0<\beta<\frac{b-a}{b+L}<1$ according to {\bf Assumption  (iv)}. In both cases, we get (\ref{relative error condition}) by noting that the existence of a positive lower bound $t^*$ of  $\{t_n\}$ is guaranteed by Proposition \ref{positive bound}.     This completes the proof. \hfill $\Box$

\begin{theorem}\label{convergence} Assume that  the  sequence  $\{x_n\}$  generated by the Algorithm satisfies  $f_{ab}(x_n)>0$ for all $n$. Under {\bf Assumptions (i)-(iv)}, the following assertions hold:
\begin{description}
  \item[{\bf (a)}] The sequence $ x_n $ has a finite length, i.e.,  $\sum_{n=0}^{+\infty} ||x_{n+1}-x_n||<+\infty$.
    \item[{\bf (b)}] The sequence  $f_{ab}(x_n)$   converges {\rm Q}-linearly to 0.
  \item[{\bf (c)}] The sequence $x_n $ converges {\rm R}-linearly to a solution $\bar{x}$ of (VIP).
\end{description}
\end{theorem}
\noindent{\bf Proof.}  From Proposition \ref{prop-key conditions}, it follows that (\ref{sufficient decrease condition}) and (\ref{relative error condition}) holds with $M_1:=\min\{\tau, b-a-\alpha\}$,  $M_2:=L+b+\frac{b-a}{\beta}$ and  $t^*$ being a  positive lower bound of  $\{t_n\}$.
By  {\bf Assumption  (i)},  the level set $[f_{ab}\leq f_{ab}(x_0)]$ is bounded, which,  together with (\ref{sufficient decrease condition}), implies that
 the sequence $\{x_n\}$ is also  bounded. Denote by $\bar{x}$ any cluster  point of the sequence $\{x_n\}$. By {\bf Assumption  (iii)}, $f$ satisfies the KL inequality at $\bar{x}$ with an exponent of $\frac{1}{2}$.   In view of these facts and the continuity of $f_{ab}$, we confirm that the sequence $\{x_n\}$  satisfies the assumptions {\bf (H1)} and {\bf (H3)}   and a variant of the assumption {\bf(H2)} in \cite{bot1}. Note that  the assumption {\bf (H2)} in \cite{bot1} requires that  $d(0, \partial f_{ab}(x_{n+1}))$,  instead of $d(0, \partial f_{ab}(x_n))$, has  an upper estimate as in the form of (\ref{relative error condition}).  In this case,  \cite[Theorem 2.9]{bot1} cannot be applied directly,  but we  can still follow the proof of \cite[Theorem 2.9]{bot1} to deduce the following: (i) {\bf (a)} holds; (ii) $x_n\to \bar{x}$ and $f_{ab}(x_n)\to f_{ab}(\bar{x})$ as $n$ goes to $\infty$; and (iii) $0\in \partial f_{ab}(\bar{x})$. In view of {\bf Assumption  (iii)} and Lemma \ref{lem-suff}, we have $\pi_a(\bar{x})=\pi_b(\bar{x})$. Then by Corollary \ref{coro-ab} {\bf (b)}, $\bar{x}$ is a solution of (VIP) or equivalently  $f_{ab}(\bar{x})=0$ (cf. Lemma \ref{lem-basic} {\bf (c)}).

It remains to show the convergence rate.
By the line search strategy in Step 3 of the Algorithm,   the following hold for all $n$:
 \begin{equation}\label{jlbswt}
 \|d_n\|\geq \beta\|x_n-\pi_a(x_n)\|,
 \end{equation}
 and
 \begin{equation}\label{dandiaoxiaj}
 \begin{array}{lll}
f_{ab}(x_{n+1})-f_{ab}(x_n)&\leq &-\min\{\tau, b-a-\alpha\} t_n\|d_n\|^2\\[0.1cm]
&\leq & -\min\{\tau, b-a-\alpha\} t^*\|d_n\|^2 \\[0.1cm]
&<&0.
\end{array}
\end{equation}
In view of (\ref{jlbswt}), we get from  Lemma \ref{lem-basic} {\bf (a)} that
$
\|d_n\|^2\geq \frac{2\beta^2}{b-a} f_{ab}(x_n)$,
which, together with (\ref{dandiaoxiaj}) and the definition of $M_1$, implies that
\[
f_{ab}(x_{n+1}) \leq  -M_1 t^* ||d_n||^2+f_{ab}(x_n)\leq (1-\frac{2\beta^2 M_1 t^*}{b-a}) f_{ab}(x_n),
\]
and hence that,
\begin{equation}\label{q-rate}
\frac{f_{ab}(x_{n+1})}{f_{ab}(x_{n})}\leq 1-\frac{2\beta^2 M_1 t^*}{b-a}=:\eta.
\end{equation}
Clearly, we have $0<\eta<1$. Then by definition \cite[pp.619-620]{nw06}, the sequence $f_{ab}(x_n)$ converges Q-linearly to 0. That is,  {\bf (b)} follows.

By the triangle inequality, the following holds for all positive integers $n$ and $m$ with $m>n$: $
\|x_n-\bar{x}\| \leq \sum_{k=n}^m\|x_{k+1}-x_k\|+\|x_{m+1}-\bar{x}\|$.
In view of {\bf (a)} and the fact that $\|x_{m+1}-\bar{x}\|\to 0$ as $m\to \infty$, we  have $\sum_{k=n}^m\|x_{k+1}-x_k\|+\|x_{m+1}-\bar{x}\|\to  \sum_{k=n}^\infty\|x_{k+1}-x_k\|$ as  $m\to \infty$, and hence $\|x_n-\bar{x}\|\leq \sum_{k=n}^\infty\|x_{k+1}-x_k\|$. In view of (\ref{sufficient decrease condition}) and (\ref{q-rate}), we further have
$$
\|x_n-\bar{x}\|\leq \sum_{k=n}^\infty\sqrt{\frac{f_{ab}(x_k)}{M_1}}\leq \sqrt{\frac{f_{ab}(x_n)}{M_1}}\sum_{k=0}^\infty\sqrt{\eta^k}=\sqrt{\frac{f_{ab}(x_n)}{M_1}}\frac{1}{1-\sqrt{\eta}}=:\zeta_n,
$$
and
 $
\frac{\zeta_{n+1}}{\zeta_n}=\sqrt{\frac{f_{ab}(x_{n+1})}{f_{ab}(x_n)}}\leq  \sqrt{\eta}.
$
As $0<\eta<1$, we have $0<\sqrt{\eta}<1$. Then by definition \cite[pp.619-620]{nw06}, $\zeta_n$ converges Q-linearly  to 0, and   $x_n$ converges R-linearly to $\bar{x}$.
  This completes the proof.
  \hfill $\Box$


\end{document}